\documentclass[a4paper,12pt]{amsart}

\usepackage{amssymb,enumerate,psfrag,graphicx,amsfonts,amsrefs,amsthm,mathrsfs,amsmath,amscd,version,graphicx}
\usepackage{xcolor}
\usepackage{tikz-cd}
\usepackage{tikz}
\usetikzlibrary{arrows}
\tikzset{
    vertex/.style={draw,circle,inner sep=2 pt, minimum size=6pt},
    edge/.style={thick},
    dedge/.style ={->,> = latex',thick}
    }
\usetikzlibrary{decorations.markings}
\usetikzlibrary{arrows.meta}
\definecolor{mygray}{gray}{0.8}

\DeclareMathOperator{\Spec}{Spec}

\newtheorem{theo}{Theorem}[section]
\newtheorem{theorem}[theo]{Theorem}
\newtheorem{corollary}[theo]{Corollary}
\newtheorem{lemma}[theo]{Lemma}
\newtheorem{proposition}[theo]{Proposition}

\theoremstyle{definition}

\renewcommand\ell{l}

\newtheorem{example}[theo]{Example}
 \newtheorem{remark}[theo]{Remark}

 \newtheorem{definition}[theo]{Definition}

\numberwithin{equation}{section}

\begin{document}
\keywords{Line graph, Signed graph, Gain graph, Spectrum of a signed graph, Forbidden subgraph, Switching isomorphism.}

\title{Characterizations of line graphs in signed and gain graphs}

\author[M. Cavaleri]{Matteo Cavaleri}
\address{Matteo Cavaleri, Universit\`{a} degli Studi Niccol\`{o} Cusano - Via Don Carlo Gnocchi, 3 00166 Roma, Italia}
\email{matteo.cavaleri@unicusano.it}

\author[D. D'Angeli]{Daniele D'Angeli}
\address{Daniele D'Angeli, Universit\`{a} degli Studi Niccol\`{o} Cusano - Via Don Carlo Gnocchi, 3 00166 Roma, Italia}
\email{daniele.dangeli@unicusano.it}

\author[A. Donno]{Alfredo Donno}
\address{Alfredo Donno, Universit\`{a} degli Studi Niccol\`{o} Cusano - Via Don Carlo Gnocchi, 3 00166 Roma, Italia}
\email{alfredo.donno@unicusano.it   (Corresponding Author)}

\begin{abstract}
We generalize three classical characterizations of line graphs to  line graphs of  signed and gain graphs: the Krausz's characterization, the van Rooij and Wilf's characterization and the Beineke's characterization. In particular, we present a list of forbidden gain subgraphs characterizing the class of gain-line graphs. In the case of a signed graph whose underlying graph is a line graph, this list consists of exactly four signed graphs. Under the same hypothesis, we prove that a signed graph is the line graph of a signed graph if and only if its eigenvalues are either greater than $-2$, or less than $2$, depending on which particular definition of line graph is adopted.
\end{abstract}

\maketitle

\begin{center}
{\footnotesize{\bf Mathematics Subject Classification (2010)}: 05C22, 05C25, 05C50, 05C76.}
\end{center}

\section{Introduction}
This article  aims at characterizing those gain graphs that are \emph{line graphs of gain graphs}. We do it by providing several equivalent conditions, inspired from the classical theory \cite{beineke, Krausz, wilf}: in terms of the existence of a partition of the edge set into (anti)balanced cliques; in terms of the $K_{1,3}$-freeness, with some supplementary conditions on the induced  triangles; by giving a list of \emph{forbidden induced gain subgraphs}. In the particular case of signed graphs, we are also able to provide a spectral characterization of such graphs.

Signed graphs were introduced in \cite{Harary}. Roughly speaking, they are graphs whose edges can be positive or negative, and interest in them goes beyond graph theory, since they can be a model for a system of interactions that can be positive or negative. A signed graph is a pair $(\Gamma,\sigma)$ where $\Gamma=(V_\Gamma,E_\Gamma)$ is the underlying graph and $\sigma\colon E_\Gamma\to\mathbb T_2=\{\pm1\}$ is the \emph{signature}.  The spectrum of the signed graph $(\Gamma,\sigma)$ is, by definition, the spectrum of its adjacency matrix \cite{zasmat}.
There is a natural \emph{switching} action on the signatures of a graph \cite{zasgraph}, inspired from Seidel's switching \cite{sei}: this operation switches  the sign of each edge with exactly one endpoint in a fixed subset of vertices. This operation induces an equivalence relation, called switching equivalence, on the set of signatures of a given graph, and it preserves the spectrum. The composition of the switching equivalence with a graph isomorphism leads to the notion of switching isomorphism.

A natural generalization of a signed graph, from the group $\mathbb T_2$ to any group $G$, is the \emph{$G$-gain graph} $(\Gamma,\psi)$, or gain graph over $G$. Here $\psi$ is a map, called the \emph{gain function}, assigning a group element to each orientation of each edge of $\Gamma$, in such a way that an element and its inverse are assigned to two opposite orientations. It is worth mentioning that gain graphs can be considered as a particular case of biased graphs \cite{zaslavsky1}, and they are also strictly related to  \emph{voltage graphs}, which are a largely investigated topic in topological graph theory \cite{voltageottobre}. The reader is referred to \cite{zasglos,zasbib} for a rich and periodically updated glossary and bibliography on signed and gain graphs.\\
\indent Switching equivalence and switching isomorphism are still defined (with the right adjustments) for gain functions and gain graphs, respectively. The adjacency matrix and the spectrum are still well defined when $G$ is a subgroup of $\mathbb T$, which is the \emph{complex unit group} (see \cite{reff1}). In the general case, there is no  canonical way to define them: one possibility is to use a \emph{represented adjacency matrix} and its spectrum  \cite{nostro}.

While for (unsigned) graphs there is a scientific consensus on what a \emph{line graph} is, even though it appeared also with different names in its nearly century-long history, this is not the case for the line graph of a signed or a gain graph. In order to narrow the field, we present some
required properties for a generalization of the line graph to signed and gain graphs.
 \\A first requirement is that the underlying graph of a suitable line graph of a signed (gain) graph, must be the line graph of the underlying graph of that signed (gain) graph.
In particular, in order to define a line graph of the gain graph $(\Gamma,\psi)$, one has to construct a gain function on the (unsigned) line graph of $\Gamma$, denoted by $L(\Gamma)$.\\
Many crucial problems when dealing with gain graphs involve properties that are invariant under switching isomorphism (e.g., \emph{balance}, spectral properties, etc.). Therefore, a second requirement is that switching-isomorphic graphs must have switching-isomorphic line graphs.\\
\begin{figure}
\begin{center}
 \begin{tikzpicture}
\draw (0,0) [black,fill=black] circle (0.08 cm);
\draw (-1.1,0.8) [black,fill=black] circle (0.08 cm);
\draw (1.1,0.8) [black,fill=black] circle (0.08 cm);
\draw (0,-1) [black,fill=black] circle (0.08 cm);

\draw [-, black] (0,0) -- (-1.1,0.8);
\draw [-, black] (0,0) -- (1.1,0.8);
\draw [-, black] (0,0) -- (0,-1);

\node [align=center] at (-1,-1.2)
{$K_{1,3}$};

\draw (-1-4.4,-1) [black,fill=black] circle (0.08 cm);
\draw (1-4,-1) [black,fill=black] circle (0.08 cm);
\draw(0-4.2,0.8)[black,fill=black] circle (0.08 cm);

\draw [-, black] (-1-4.4,-1) -- (1-4,-1);
\draw [-, black] (0-4.2,0.8) -- (1-4,-1);
\draw [-, black] (-1-4.4,-1) -- (0-4.2,0.8);

\node [align=center] at (-2.2-4,-1.2)
{$K_3$};

\end{tikzpicture}
\end{center}\caption{The graphs $K_3$ and $K_{1,3}$.}\label{fig:KK}
\end{figure}
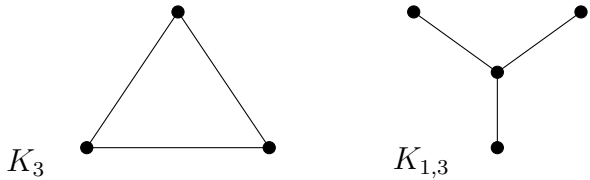
\noindent Finally, a fundamental result on line graphs is the Whitney isomorphism theorem \cite{Whitney}: with the exception of the pair of graphs consisting of the complete graph $K_3$ and of the complete bipartite graph $K_{1,3}$ (see Fig. \ref{fig:KK}), two graphs are isomorphic if and only if their line graphs are isomorphic. In particular, if a graph is not $K_3$ and it is a line graph, it is possible to reconstruct, up to isomorphism, the graph for which it is the line graph of. If one requires the analogous possibility for the line graph of a gain graph, the third requirement is that only switching equivalent gain functions induce line graphs with switching equivalent gain functions. \\
There exist some definitions of line graphs for signed graphs in the literature (e.g. \cite{ajk, beh}), which do not satisfy one or more of the aforementioned conditions.

To the best of our knowledge, the first definition with the required properties is given in \cite{zas84} for signed graphs. The same definition is considered in \cite{ger, zasmat} and it is generalized in
\cite{reff0} to gain graphs over an Abelian group. Another suitable definition for the line graph of a signed graph appears in \cite{francescolap}, and it is considered also in \cite{francescomp,francescostar} and generalized to gain graphs over the complex unit group $\mathbb{T}$ \cite{fm,francesco}.
These definitions differ only by a switch of sign (see \cite{totalarxiv} for remarks on this regard) or, more generally, by a multiplication of the gain functions by a \emph{central involution $s$} of the group $G$.
By introducing two central involutions $s_1$ and $s_2$ of $G$ as parameters, \cite{line} gives a unifying generalization of the line graph construction to gain graphs over an arbitrary group $G$, providing an answer to a question formulated by N. Reff in \cite{reff0} about a possible extension of this construction to the non-Abelian case.
All these definitions pass through the choice of a $G$-phase ($G$-incidence phase function in the language of \cite{reff0}). Actually, the line graph is more properly defined for
switching equivalence classes: in \cite{line} this is implemented by the map $\mathcal L$ from the switching equivalence classes of $\Gamma$ to those of $L(\Gamma)$ (see Theorem \ref{theo:vecchio} of the present paper).\\
\indent The interest in gain line graphs goes beyond the issue related to giving a suitable definition, and it is also focused on their spectral properties. Indeed, it turns out that, as in the classical setting, the adjacency matrix of a gain graph, the adjacency matrix of its gain line graph, and the $G$-phase matrix are closely related \cite{francesco, line, reff0}. In the very recent paper \cite{belardoquaternions}, a spectral investigation of quaternion unit gain graphs and associated line graphs has been developed.

\begin{figure}
\begin{center}
 \begin{tikzpicture}
\draw (0,0) [black,fill=black] circle (0.08 cm);
\draw (-1,0.7) [black,fill=black] circle (0.08 cm);
\draw (1,0.7) [black,fill=black] circle (0.08 cm);
\draw (0,-1) [black,fill=black] circle (0.08 cm);

\draw [-, black] (0,0) -- (-1,0.7);
\draw [-, black] (0,0) -- (1,0.7);
\draw [-, black] (0,0) -- (0,-1);

\node [align=center] at (-1,-1)
{$G_1$};

\draw (-1+3.5,0) [black,fill=black] circle (0.08 cm);
\draw (1+3.5,0) [black,fill=black] circle (0.08 cm);
\draw (0+3.5,1) [black,fill=black] circle (0.08 cm);
\draw (0+3.5,-1) [black,fill=black] circle (0.08 cm);
\draw (1.5+3.5,0) [black,fill=black] circle (0.08 cm);

\draw [-, black] (-1+3.5,0) -- (1+3.5,0);
\draw [-, black] (-1+3.5,0) -- (0+3.5,1);
\draw [-, black] (-1+3.5,0) -- (0+3.5,-1);
\draw [-, black] (1+3.5,0) -- (0+3.5,1);
\draw [-, black] (1+3.5,0) -- (0+3.5,-1);
\draw [-, black] (1.5+3.5,0) -- (0+3.5,1);
\draw [-, black] (1.5+3.5,0) -- (0+3.5,-1);

\node [align=center] at (-1+3.5,-1)
{$G_2$};

\draw (-1+7.4,-1) [black,fill=black] circle (0.08 cm);
\draw (1+7,-1) [black,fill=black] circle (0.08 cm);
\draw(0+7.2,1)[black,fill=black] circle (0.08 cm);

\draw(0+7.2,0.4)[black,fill=black] circle (0.08 cm);
\draw(0+7.2,-0.2)[black,fill=black] circle (0.08 cm);

\draw [-, black] (-1+7.4,-1) -- (1+7,-1);
\draw [-, black] (0+7.2,1) -- (1+7,-1);
\draw [-, black] (-1+7.4,-1) -- (0+7.2,1);

\draw [-, black] (0+7.2,0.4) -- (1+7,-1);
\draw [-, black] (-1+7.4,-1) -- (0+7.2,0.4);

\draw [-, black] (0+7.2,-0.2) -- (1+7,-1);
\draw [-, black] (-1+7.4,-1) -- (0+7.2,-0.2);

\draw [-, black] (0+7.2,-0.2) -- (0+7.2,0.4);

\draw [-, black] (0+7.2,1) -- (0+7.2,0.4);

\node [align=center] at (-1+7,-1)
{$G_3$};

\draw (-1+0,0-3) [black,fill=black] circle (0.08 cm);
\draw (1+0,0-3) [black,fill=black] circle (0.08 cm);
\draw (0+0,1-3) [black,fill=black] circle (0.08 cm);
\draw (0+0,-1-3) [black,fill=black] circle (0.08 cm);

\draw (1.5+0,1-3) [black,fill=black] circle (0.08 cm);
\draw (1.5+0,-1-3) [black,fill=black] circle (0.08 cm);

\draw [-, black] (-1+0,0-3) -- (1+0,0-3);
\draw [-, black] (-1+0,0-3) -- (0+0,1-3);
\draw [-, black] (-1+0,0-3) -- (0+0,-1-3);
\draw [-, black] (1+0,0-3) -- (0+0,1-3);
\draw [-, black] (1+0,0-3) -- (0+0,-1-3);
\draw [-, black] (1.5+0,1-3) -- (0+0,1-3);
\draw [-, black] (1.5+0,-1-3) -- (0+0,-1-3);

\node [align=center] at (-1+0,-1-3)
{$G_4$};

\draw (-1+0,0-6) [black,fill=black] circle (0.08 cm);
\draw (1+0,0-6) [black,fill=black] circle (0.08 cm);
\draw (0+0,1-6) [black,fill=black] circle (0.08 cm);
\draw (0+0,-1-6) [black,fill=black] circle (0.08 cm);

\draw (1.5+0,1-6) [black,fill=black] circle (0.08 cm);
\draw (1.5+0,-1-6) [black,fill=black] circle (0.08 cm);

\draw [-, black] (-1+0,0-6) -- (1+0,0-6);
\draw [-, black] (-1+0,0-6) -- (0+0,1-6);
\draw [-, black] (-1+0,0-6) -- (0+0,-1-6);
\draw [-, black] (1+0,0-6) -- (0+0,1-6);
\draw [-, black] (1+0,0-6) -- (0+0,-1-6);
\draw [-, black] (1.5+0,1-6) -- (0+0,1-6);
\draw [-, black] (1.5+0,-1-6) -- (0+0,-1-6);
\draw [-, black] (1.5+0,-1-6) --  (1.5+0,1-6);

\node [align=center] at (-1+0,-1-6)
{$G_7$};

\draw (-1+3.5,0-3) [black,fill=black] circle (0.08 cm);
\draw (1+3.5,0-3) [black,fill=black] circle (0.08 cm);
\draw (0+3.5,1-3) [black,fill=black] circle (0.08 cm);
\draw (0+3.5,-1-3) [black,fill=black] circle (0.08 cm);

\draw (1.5+3.5,-1-3) [black,fill=black] circle (0.08 cm);
\draw (0+3.5,0+0.4-3) [black,fill=black] circle (0.08 cm);

\draw [-, black] (-1+3.5,0-3) -- (1+3.5,0-3);
\draw [-, black] (-1+3.5,0-3) -- (0+3.5,1-3);
\draw [-, black] (-1+3.5,0-3) -- (0+3.5,-1-3);
\draw [-, black] (1+3.5,0-3) -- (0+3.5,1-3);
\draw [-, black] (1+3.5,0-3) -- (0+3.5,-1-3);

\draw [-, black] (0+3.5,0+0.4-3) -- (1+3.5,0-3);
\draw [-, black] (-1+3.5,0-3) -- (0+3.5,0+0.4-3);
\draw [-, black] (0+3.5,0+0.4-3) -- (0+3.5,1-3);

\draw [-, black] (1.5+3.5,-1-3) -- (0+3.5,-1-3);

\node [align=center] at (-1+3.5,-1-3)
{$G_5$};

\draw (-1+7,0-3) [black,fill=black] circle (0.08 cm);
\draw (1+7,0-3) [black,fill=black] circle (0.08 cm);
\draw (0+7,1-3) [black,fill=black] circle (0.08 cm);
\draw (0+7,-1-3) [black,fill=black] circle (0.08 cm);

\draw (0+7,0+0.4-3) [black,fill=black] circle (0.08 cm);
\draw (0+7,0-0.4-3) [black,fill=black] circle (0.08 cm);

\draw [-, black] (-1+7,0-3) -- (1+7,0-3);
\draw [-, black] (-1+7,0-3) -- (0+7,1-3);
\draw [-, black] (-1+7,0-3) -- (0+7,-1-3);
\draw [-, black] (1+7,0-3) -- (0+7,1-3);
\draw [-, black] (1+7,0-3) -- (0+7,-1-3);

\draw [-, black] (0+7,0+0.4-3) -- (1+7,0-3);
\draw [-, black] (-1+7,0-3) -- (0+7,0+0.4-3);
\draw [-, black] (0+7,0+0.4-3) -- (0+7,1-3);

\draw [-, black] (0+7,0-0.4-3) -- (1+7,0-3);
\draw [-, black] (-1+7,0-3) -- (0+7,0-0.4-3);
\draw [-, black] (0+7,0-0.4-3) -- (0+7,-1-3);

\node [align=center] at (-1+7,-1-3)
{$G_6$};

\draw (-0.5+3.5,1-6) [black,fill=black] circle (0.08 cm);
\draw (-0.5+3.5,0-6) [black,fill=black] circle (0.08 cm);
\draw (-0.5+3.5,-1-6) [black,fill=black] circle (0.08 cm);
\draw (1+3.5,1-6) [black,fill=black] circle (0.08 cm);
\draw (1+3.5,0-6) [black,fill=black] circle (0.08 cm);
\draw (1+3.5,-1-6) [black,fill=black] circle (0.08 cm);

%orizz
\draw [-, black] (-0.5+3.5,1-6) -- (1+3.5,1-6);
\draw [-, black] (-0.5+3.5,0-6) -- (1+3.5,0-6);
\draw [-, black] (-0.5+3.5,-1-6) -- (1+3.5,-1-6);

%verticali
\draw [-, black] (-0.5+3.5,-1-6) -- (-0.5+3.5,0-6);
\draw [-, black] (-0.5+3.5,1-6) -- (-0.5+3.5,0-6);

\draw [-, black] (1+3.5,-1-6) -- (1+3.5,0-6);
\draw [-, black] (1+3.5,1-6) -- (1+3.5,0-6);

%diagonali
\draw [-, black] (-0.5+3.5,0-6) -- (1+3.5,1-6);
\draw [-, black] (-0.5+3.5,-1-6) -- (1+3.5,0-6);

\node [align=center] at (-1+3.5,-1-6)
{$G_8$};

\draw (0+7,0-6) [black,fill=black] circle (0.08 cm);
\draw (0+7-0.5,0-6-1) [black,fill=black] circle (0.08 cm);
\draw (0+7+0.5,0-6-1) [black,fill=black] circle (0.08 cm);

\draw (0+7-1,0-6+0.2) [black,fill=black] circle (0.08 cm);
\draw (0+7+1,0-6+0.2) [black,fill=black] circle (0.08 cm);

\draw (0+7,0-6+1) [black,fill=black] circle (0.08 cm);

%stella
\draw [-, black] (0+7,0-6) -- (0+7-0.5,0-6-1);
\draw [-, black] (0+7,0-6) --(0+7+0.5,0-6-1);
\draw [-, black] (0+7,0-6) -- (0+7-1,0-6+0.2);
\draw [-, black] (0+7,0-6) --  (0+7+1,0-6+0.2);
\draw [-, black] (0+7,0-6) -- (0+7,0-6+1);

%ciclo
\draw [-, black] (0+7-0.5,0-6-1) --(0+7+0.5,0-6-1);
\draw [-, black] (0+7+1,0-6+0.2) --(0+7+0.5,0-6-1);
\draw [-, black] (0+7+1,0-6+0.2) --(0+7,0-6+1);
\draw [-, black] (0+7,0-6+1) --(0+7-1,0-6+0.2);
\draw [-, black] (0+7-0.5,0-6-1) --(0+7-1,0-6+0.2);

\node [align=center] at (-1+7,-1-6)
{$G_9$};

\end{tikzpicture}
\end{center}\caption{The list of nine forbidden subgraphs $\mathcal X$.}\label{fig:9}
\end{figure}
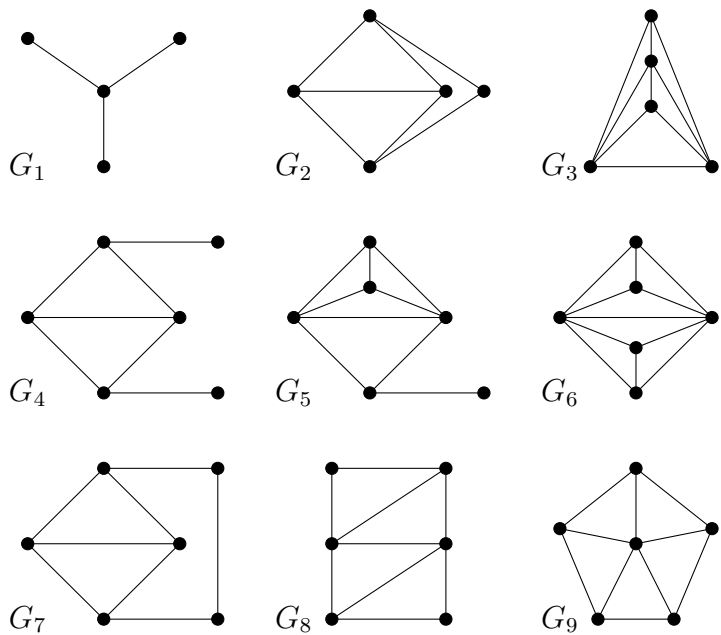

The map $\mathcal{L}$ of \cite{line} is proved to be injective, giving a \emph{generalization of Whitney isomorphism theorem} for line graphs of gain graphs (that holds in all of the previous compatible definitions). The Whitney isomorphism theorem is just the first of a series of results about line graphs and the aim of this paper is to take further steps in this generalization from the classical theory to that of gain graphs. Beineke in \cite{beineke} gives a list $\mathcal X$ of nine graphs, that we represent in Fig.~\ref{fig:9}, called \emph{forbidden subgraphs for a line graph}, with the property that a graph $\Gamma$ is the line graph of a simple graph if and only if none of its induced subgraphs is isomorphic to one of the graphs in $\mathcal X$, briefly, $\Gamma$ is \emph{$\mathcal X$-free}. This characterization is particularly useful, also from an algorithmic point of view,
since one can establish whether $\Gamma$ is a line graph or not just by looking at its induced subgraphs of at most $6$ vertices, which is the upper bound for the number of vertices of the graphs in $\mathcal X$. Analogous characterizations are given for \emph{generalized line graphs} in \cite{sin82,sin81} and for line graphs of multigraphs \cite{bermond}. The original path for the proof of Beineke's characterization in \cite{beineke} passes through Krausz's characterization \cite{Krausz}, that is given in terms of a partition of the edge set into complete subgraphs, or \emph{cliques}; and through van Rooij and Wilf's characterization, that is given in terms of the $K_{1,3}$-freeness and odd triangles (see Definition \ref{defiodd9}). These characterizations can be summarized in the following theorem (see also \cite{Bharary}).
\begin{theorem}\label{theo:0}
A connected graph $L=(V_L,E_L)$ is the line graph of a simple graph if and only if it satisfies one of the following equivalent conditions.
\begin{itemize}
\item[\underline{\footnotesize{Krausz:}}] There exists a partition of the edge set $E_L=E_1\sqcup E_2 \sqcup\cdots\sqcup E_k$ such that every vertex of $L$ is endpoint of edges from at most two elements of the partition and the subgraph $L_{E_i}$ of $L$ induced by $E_i$ is complete, for each $i=1,\ldots, k$.\\
\item[\underline{\footnotesize{van Rooij and Wilf:}}]
 $L$ is $K_{1,3}$-free and, if $T_1$, $T_2$ are adjacent odd triangles, their vertices induce a subgraph  of $L$ that is complete.\\
\item[\underline{\footnotesize{Beineke:}}] $L$ is  $\mathcal X$-free.
\end{itemize}
\end{theorem}

In Theorem \ref{theo:1} we give the generalization of Theorem \ref{theo:0} to gain graphs, i.e., we give the analogous necessary and sufficient conditions for a gain graph to be a  line graph of a gain graph. According with the definition of gain line graph given in \cite{line}, which extends those used in \cite{francesco, fm, francescolap, ger, reff0} for signed and complex unit gain graphs, we will consider only gain graphs whose underlying graph is simple.
Since the underlying graph of a gain-line graph is a line graph itself, the classical conditions remain necessary, but our new conditions in the gain graphs setting are stronger. For the analogue of Krausz's characterization, we need the further condition that the complete graphs induced by the edge partition  have induced gain functions of a specific form.
In the generalized van Rooij and Wilf condition, the additional requirement concerns  the gain of each odd triangle, together with one more condition that is not trivially satisfied when the underlying graph is isomorphic to one among the three graphs $F_1,F_2,F_3$ of Fig.~\ref{fig:tre}.
Finally,  in the generalization of Beineke's condition, we replace the list of forbidden subgraphs $\mathcal X$ with a list of \emph{forbidden gain subgraphs} $\mathcal Y$ (see Eq.~\eqref{eq:y}).\\
\indent If one restricts the problem to recognize if a gain graph $(L(\Gamma),\zeta)$ is a gain-line graph, under the assumption that the underlying graph $L(\Gamma)$ is a line graph, something more can be said (see Corollary \ref{cor:primo}).

There exist several operations on signed and gain graphs leaving the spectrum unchanged (see \cite{cos,gm}). As we already observed, the switching is one of them \cite{zasmat}. It follows that the spectrum of a switching equivalence class of signed graphs does not depend on the particular choice of a representative; therefore, the spectrum of the line graph of a signed graph is well defined.
As a consequence of  Theorem \ref{theo:1}, in the case of signed graphs, one can recognize which signatures on a line graph induce line graphs of signed graphs only by means of their spectra. This is shown in Theorem  \ref{cor:spec}, together with a characterization of signed line graphs via $4$ forbidden signed subgraphs (Fig. \ref{fig:fsegn} and Fig. \ref{fig:fsegnmeno}).
In \cite{meno2,vi87,survey}  the class of \emph{signed graphs represented by $D_\infty$} is also characterized in an analogous way. This class coincides with that of \emph{reduced line graphs} of simple signed graphs admitting   parallel edges, though only with opposite sign  \cite{zas84,zasmat}. A comparison with Theorem  \ref{cor:spec} is given in Remark \ref{rem:vi}.

Finally, we characterize in Corollary \ref{cor:cp} the connected graphs $\Gamma$ that, equipped with every gain function, are gain-line graphs. If the group $G$ is nontrivial, they are exactly the \emph{cycles}  and the \emph{paths}. As a consequence, if $\Gamma$ is not a cycle or a path, the \emph{circuit rank} of $L(\Gamma)$ is strictly greater than that of $\Gamma$ (see Corollary \ref{cor:cp1}).  Moreover, the cyclic graphs and the  path graphs are  exactly the connected line graphs  with all eigenvalues of modulus at most $2$ (see Corollary \ref{cor:cp2}).

\section{Gain graphs, $G$-phases and gain-line graphs}\label{sec:prel}
Let $\Gamma=(V_\Gamma,E_\Gamma)$ be a finite, connected, simple, undirected graph, with at least one edge. The set $V_\Gamma$ is the vertex set, and the set $E_\Gamma$ is the edge set, consisting of unordered pairs of the type $\{u,v\}$, with $u,v\in V_\Gamma
$. We write $u\sim v$ if $\{u,v\}\in E_\Gamma$, then  we say that $u$ and $v$ are \emph{adjacent} and that are \emph{endpoints} of the edge $\{u,v\}$.
We will use the set theoretic notation for the edges: for $v\in V_\Gamma$ and $e\in E_\Gamma$ we write $v\in e$ if the edge $e$ is \emph{incident} to $v$, that is, $v$ is one of the endpoints of $e$. If $e_1,e_2\in E_\Gamma$ are both incident to a vertex, we denote that vertex as $e_1\cap e_2$. We write $e_1\cap e_2=\emptyset$ if $e_1$ and $e_2$ do not share a common vertex. The \emph{line graph} $L(\Gamma)$ is the graph with vertex set $E_\Gamma$, whose vertices $e_1,e_2$ are adjacent if $e_1\cap e_2\neq \emptyset$ in $\Gamma$.

Let $G$ be a group and consider a map
$\psi\colon \{(u,v)\in V_\Gamma^{2} \mid u\sim v\}\to G$ such that $\psi(u,v)=\psi(v,u)^{-1}$. The pair $(\Gamma,\psi)$ is a \emph{$G$-gain graph}  (or equivalently, a gain graph over $G$) and $\psi$ is said to be a \emph{gain function}. The graph $\Gamma$ is the \emph{underlying graph} of the gain graph $(\Gamma,\psi)$.
We denote by $G(\Gamma)$ the set of all gain functions of $\Gamma$ over $G$. The most studied gain graphs are those over the group $\mathbb T_2=\{\pm1\}$ (more properly called \emph{signed graphs}) or over the group $\mathbb T=\{z\in \mathbb C \mid |z|=1\}$  (more properly called \emph{complex unit gain graphs}).

Let $W$ be a \emph{walk} of length $\ell$ in $\Gamma$, that is, an ordered sequence of $\ell+1$  vertices of $\Gamma$, say $v_0,v_1,\ldots, v_\ell$, with $v_i\sim v_{i+1}$. We will denote by $|W|$ the length of the walk $W$. The gain of $W$ is then defined as
$$
\psi(W):=\psi(v_0,v_1)\cdots \psi(v_{\ell-1},v_\ell).
$$
A \emph{closed walk} of length $\ell$ is a walk of length $\ell$ with $v_0=v_{\ell}$.

We denote by $1_G$ the identity element of $G$ and we say that $s\in G$  is an \emph{involution} if $s^2=1_G$. An element $g\in G$ is said to be \emph{central} if it commutes with any other element of $G$.
If $s\in G$ is an involution, we denote, in bold, by $\bold{s}$  the map such that $\bold{s}(u,v)=s$  whenever $u\sim v$. It is clear that $\bold{s}$ is a gain function for $\Gamma$.
The gain graph $(\Gamma,\psi)$ is said to be \emph{balanced} if $\psi(W)=1_G$ for every closed walk $W$. For example $(\Gamma, \bold{1_G})$ is trivially balanced.
Notice that, if $T$ is a tree, one has $\psi(W)=1_G$ for every $\psi\in G(T)$ and any closed walk $W$ of $T$, hence a gain tree $(T,\psi)$ is always balanced.
\begin{definition}\label{def:induG}
Let $(\Gamma,\psi)$ be a gain graph, with $\Gamma=(V_\Gamma, E_\Gamma)$. Let $A\subseteq V_\Gamma$. Then:
\begin{itemize}
\item
$E_\Gamma(A)\subseteq E_\Gamma$ is the subset of edges whose endpoints are both in $A$;
\item
$\Gamma_A=(A,E_\Gamma(A))$ is the \emph{subgraph induced} by the subset $A$ of vertices;
\item
$\psi_{\Gamma_A}\in G(\Gamma_A)$, or simply $\psi_A\in G(\Gamma_A)$, is the gain function on $\Gamma_A$ such that $\psi_A(u,v)=\psi(u,v)$ for every $u,v\in A$ with $u \sim v$.
\end{itemize}
Let $B\subseteq E_\Gamma$. Then:
\begin{itemize}
\item  $V_\Gamma(B)\subseteq V_\Gamma$ is the subset of vertices that are endpoints of edges in $B$;
\item $\Gamma_B=(V_\Gamma(B),B)$ is the \emph{edge-induced subgraph} by the subset $B$ of edges;
\item $\psi_{\Gamma_B}\in G(\Gamma_B)$, or simply $\psi_B\in G(\Gamma_B)$,  is the gain function  on $\Gamma_B$ such that $\psi_B(u,v)=\psi(u,v)$ when $\{u,v\}\in B$.
\end{itemize}
\end{definition}
\begin{remark}\label{rem:ab}
For $B\subseteq E_\Gamma$  and $A:=V_\Gamma(B)\subseteq V_\Gamma$,
the graph $\Gamma_B$ differs from $\Gamma_A$ in general. This happens, for instance, if $\Gamma$ is not a tree, and $B$ is a spanning tree of $\Gamma$. However, if $B$ is also equal to $E_\Gamma(A)$, then $\Gamma_A=\Gamma_B$. This is the case, for example, when $\Gamma_B$ is a complete graph.
\end{remark}

A \emph{cycle} $C$ in $\Gamma$ is a subgraph induced by $k$ of its edges that is isomorphic to the cyclic graph $C_k$.
A cycle involving $3$ edges is called \emph{triangle}. As noticed in the previous remark, a triangle $T$ can be also seen as a subgraph induced by $3$ adjacent vertices $v_1,v_2,v_3$. Two distinct triangles $T_1,T_2$ of $\Gamma$ are said to be \emph{adjacent} if they share an edge.
\begin{definition}\cite{beineke}   \label{defiodd9}
A triangle $T$ of $\Gamma$ on the vertices $v_1,v_2,v_3$ is said to be an \emph{even triangle} if the number $|\{v\in\{v_1,v_2,v_3\} \mid v\sim w\}|$ is even for every $w\in V_\Gamma$; it is said to be an \emph{odd triangle} otherwise.
\end{definition}
In other words, a triangle $T$ is odd if there exists at least one vertex $w\in V_\Gamma$ which is adjacent to an odd number of vertices of $T$.

A fundamental concept in the theory of gain graphs, inherited from the theory of signed graphs, is the \emph{switching equivalence} and the  \emph{switching isomorphism}.

\begin{definition}\label{defswe}
Two gain functions $\psi_1$ and $\psi_2$ on the same underlying graph
$\Gamma =(V_\Gamma, E_\Gamma)$ are switching equivalent, and we  write $\psi_1\sim \psi_2$, if there exists a \emph{switching function} from $\psi_1$ to $\psi_2$, that is,
 $f\colon V_\Gamma\to G$ such that
\begin{equation}\label{eqsw}
\psi_2(u,v)=f(u)^{-1} \psi_1(u,v)f(v), \qquad \forall u,v\in V_\Gamma, u\sim v.
\end{equation}
Two gain graphs $(\Gamma_1,\psi_1)$ and $(\Gamma_2,\psi_2)$ are \emph{switching isomorphic} if there is a graph-isomorphism $\phi\colon V_{\Gamma_1}\to V_{\Gamma_2}$ such that
$\psi_1\sim (\psi_2\circ \phi)$, where $\psi_2\circ \phi$ is the gain function on $\Gamma_1$ such that $(\psi_2\circ \phi)(u,v)=\psi_2(\phi(u),\phi(v))$.
\end{definition}
When Eq. \eqref{eqsw} holds we shortly write $\psi_2=\psi_1^{f}$. We denote by $[\psi]$ the \emph{switching equivalence class} of the gain function $\psi$ and by $[G(\Gamma)]$ the set of all switching equivalence classes of gain functions on $\Gamma$.

As a consequence of Eq. \eqref{eqsw}, if $\psi_1\sim \psi_2$ and $W$ is a closed walk, then $\psi_1(W)$ and $\psi_2(W)$ are conjugated elements in $G$ (see \cite[Proposition 2.1]{reff0}).
Moreover, a gain graph $(\Gamma, \psi)$ is balanced if and only if $ \psi\sim  \bold{1_G}$ (see \cite[Lemma~5.3]{zaslavsky1} or \cite[Lemma~2.1]{reff1}).
The next lemma gives a generalization  of this result.

\begin{lemma}\label{lem:s}
Let $(\Gamma,\psi)$ be a $G$-gain graph and let $s$ be a central involution of $G$. Then $\psi\sim \bold{s}$ if and only if,
for every closed walk $W$, one has $\psi(W)=s^{|W|}$.
\end{lemma}
\begin{proof}
For every closed walk $W$ we have, by definition of $\bold{s}$, that $\bold{s}(W)=s^{|W|}$. Suppose that $\psi\sim \bold{s}$, then $\psi(W)$ is conjugated to $s^{|W|}$. Since $s$ is central, we have $\psi(W)=s^{|W|}$.\\
\indent Vice versa, suppose that $\psi(W)=s^{|W|}$ for every closed walk $W$.
Choose a subset $B\subseteq E_{\Gamma}$ such that the subgraph $T:=\Gamma_B$ induced by $B$ is a spanning tree of $\Gamma$, so that $V_T =V_\Gamma$.
 Since $T$ is a tree, all of its gain functions are balanced and then $\psi_T\sim \bold{1_G}$, but also the gain function $\bold{s}$ restricted on $T$, is switching equivalent with  $\bold{1_G}$ and then by transitivity $\psi_T\sim \bold{s}$.
 Let $f\colon V_\Gamma\to G$ be the switching function such that $\psi_T^f=\bold{s}$. Consider now $\psi'\in G(\Gamma)$ with $\psi':=\psi^f$. By definition of $\psi'$ we have $\psi\sim \psi'$.
Then we are done if we prove that $\psi'=\bold{s}$.
By definition of $f$, we have that  if $\{u,v\}\in B$ then $\psi'(u,v)=\psi^f(u,v)=\psi_T^f(u,v)=s$. If $\{u,v\}\notin B$, as $T$ is a spanning tree, there is a cycle, and so a closed walk $W$, containing only edges of $B$ except for $\{u,v\}$ (suppose, in the order $u$, $v$).
As a consequence, we have that $\psi'(W)=s^{|W|-1} \psi'(u,v)$. Since $\psi$ and $\psi'$ are switching equivalent and $W$ is closed, then $\psi(W)$ and $\psi'(W)$ are conjugated.
Since $\psi(W)=s^{|W|}$ is central, it follows that  $\psi'(u,v)=s$.
The thesis follows since $\psi'=\bold{s}$ and then $\psi\sim \bold{s}$.
\end{proof}
We are going to recall the definition of line graph of a gain graph from \cite{line}, which extends the definition for the Abelian case given in \cite{reff0}.
For a graph $\Gamma$ with $n$ vertices and $m$ edges, as a generalization of the incidence matrices, we consider the space $\mathcal H_\Gamma$ of the $G$-phases. A $G$-phase $H\in \mathcal H_\Gamma$ is a map $H\colon V_\Gamma\times E_\Gamma\to G\cup \{0\}$ with the property that, for $v\in V_\Gamma$ and $e\in E_\Gamma$, we have $v\notin e$ if and only if  $H(v,e)=0$. By fixing an ordering on $V_\Gamma=\{v_1,\ldots,v_n\}$  and $E_\Gamma=\{e_1,\ldots,e_m\}$, we will interpret $H$ as an $n\times m$ matrix whose entry $H_{i,j}$ is in $G$ if $v_i\in e_j$, and $0$ otherwise. Notice that the positions of the zeros in $H$ coincide with the position of the zeros in the classical incidence matrix of $\Gamma$. Formally, these objects are particular cases of group algebra valued matrices $M_{n\times m}(\mathbb C G)$, which have been already applied to gain graphs in  \cite{nostro,line}.

In order to construct a line graph consistent with the several definitions in literature on signed, complex unit, and Abelian gain graphs, we need the freedom to choose  two parameters: two central involutions (possibly trivial) $s_1$ and $s_2$ of $G$.
The involution $s_1$ plays a crucial role in the construction of a gain function on $\Gamma$ starting from a $G$-phase $H$; the involution $s_2$ plays a crucial role in the construction of a gain function on $L(\Gamma)$.
\begin{definition}\label{def:psi}
Let $s_1$ and $s_2$ be two fixed central involutions of the group $G$. Let $\Gamma$ be a connected graph with $V_\Gamma=\{v_1,\ldots,v_n\}$  and $E_\Gamma=\{e_1,\ldots,e_m\}$. We have two maps:
\begin{equation*}
\begin{split}
&\Psi\colon \mathcal H_\Gamma \to G(\Gamma)\qquad \qquad \qquad \qquad\, \Psi_L\colon  \mathcal H_\Gamma \to G(L(\Gamma))\\
&\Psi(H)(v_i,v_j)=s_1 H_{i,k} (H_{j,k})^{-1} \qquad \Psi_L(H)(e_p,e_q)=s_2 (H_{r,p})^{-1} H_{r,q}
\end{split}
\end{equation*}
where $e_k=\{v_i,v_j\}$ and  $v_r= e_p \cap e_q.$
\end{definition}
The following theorem is a reformulation of the results of \cite[Theorem 4.25, Corollary 4.26]{line}.
\begin{theorem}\label{theo:vecchio}
For every graph $\Gamma$ there exists an injective map $\mathcal L\colon [G(\Gamma)]\to [G(L(\Gamma))]$ such that
for $\psi\in G(\Gamma)$ and $\zeta\in G(L(\Gamma))$ one has $\mathcal L([\psi])=[\zeta]$ if and only if there exists $H\in \mathcal H_\Gamma$ such that $\Psi(H)=\psi$ and $\Psi_L(H)=\zeta$.
\end{theorem}
\begin{definition}\label{defimouse}
If $\psi\in G(\Gamma)$ and $\zeta\in G(L(\Gamma))$ are such that $\mathcal L([\psi])=[\zeta]$, we say that
$(L(\Gamma),\zeta)$ is a \emph{line graph of the gain graph $(\Gamma,\psi)$}. A gain graph $(L,\zeta)$ is a \emph{gain-line graph} if it is a line graph of some gain graph.
\end{definition}
According with Definition \ref{defimouse}, the switching isomorphism class of $(L(\Gamma),\zeta)$ is the line graph of the switching isomorphism class of $(\Gamma,\psi)$.\\
\indent By virtue of Theorem \ref{theo:vecchio}, one can also state that $(L(\Gamma),\zeta)$ is a gain-line graph if and only if there exists a graph $\Gamma$ and $H\in \mathcal H_\Gamma$ such that $\Psi_L(H)=\zeta$.

\begin{example}\rm
Starting from a gain graph $(\Gamma,\psi)$ and an ordering $\{v_1,\ldots, v_n\}$ of $V_\Gamma$, we can define a particular $G$-phase $H_< \in \mathcal{H}_\Gamma$ as
\begin{equation*}
(H_{<})_{i,k}:=
\begin{cases}
0 & \mbox{ if } v_i\not\in e_k\\
\psi(v_i,v_j) &\mbox{ if } e_k=\{v_i,v_j\} \mbox { and } i<j\\
s_1 &\mbox{ if } e_k=\{v_i,v_j\} \mbox { and } i>j.
\end{cases}
\end{equation*}
By using Definition \ref{def:psi}, it is not difficult to check that $\Psi(H_<)=\psi$ and that $\Psi_L(H_<)=\zeta$, with
\begin{equation}\label{eq:line}
\zeta(\{v_j,v_i\},\{v_i,v_k\}):=
\begin{cases}
s_1s_2\psi(v_i,v_k)\quad &\mbox{ if } j<i, \,i<k\\
s_2 \quad &\mbox{ if } j<i,\, i>k\\
s_2 \psi(v_i,v_j)\psi(v_i,v_k)   \quad &\mbox{ if } j>i,\, i<k.
\end{cases}
\end{equation}
It follows that $(L(\Gamma),\zeta)$ is a line graph of the gain graph $(\Gamma,\psi)$.
For example, in the case of a signed graph with $s_1=s_2=-1$ (consistently with \cite{zasmat}) we have that the sign of an edge joining
$\{v_1,v_2\}$ and $\{v_2,v_3\}$ is the sign of $\{v_2,v_3\}$; the sign of an edge joining $\{v_2,v_3\}$ and $\{v_3,v_1\}$ is negative and finally the sign of an edge joining $\{v_3,v_1\}$ and $\{v_1,v_2\}$ is the the opposite of the  product of the sign of $\{v_3,v_1\}$ and $\{v_1,v_2\}$. Notice that, by Theorem \ref{theo:vecchio}, changing the ordering does not change the switching equivalence class of $\zeta$.
\end{example}

\begin{remark}\label{rem:s2}
If $(L,\zeta)$ is a gain-line graph, then $L$ is a line graph in the classical sense. On the other hand, a gain graph $(L,\zeta)$ with underlying graph that is a line graph
is not necessarily a gain-line graph. In other words, the map $\mathcal L$, or equivalently $\Psi_L$, is not surjective in general.
By Definition \ref{def:psi}, the range of $\Psi_L$ always contains the gain function $\bold{s_2}$ as  image of the $G$-phase with entries in $\{0,1_G\}$. As a consequence, we have  that
 $(L,\bold{s_2})$ is a gain-line graph if and only if $L$ is a line graph. In this sense, the property of being a gain-line for gain graphs is a generalization of the property of being a line graph in the classical setting.
\end{remark}
 \indent The next lemma shows that, for the line graph of a tree the class $[\bold{s_2}]$ is the only switching equivalence class of gain functions giving a gain-line graph.
\begin{lemma}\label{lem:tree}
Let $T$ be a tree and let $L(T)$ be its line graph. Then $(L(T),\zeta)$ is a gain-line if and only if $\zeta\sim \bold{s_2}$.
\end{lemma}
\begin{proof}
For any graph $\Gamma$ one can prove, for example by using Eq. \eqref{eq:line}, that $(L(\Gamma),\bold{s_2})$ is a line graph of $(\Gamma,\bold{s_1})$.
It follows that $\mathcal L([\bold{s_1}])=[\bold{s_2}]$. The thesis follows by noticing that all gain functions on a tree  are balanced and then switching equivalent with each other; in particular, there exists a unique class in $[G(T)]$.
\end{proof}

We have already defined induced graphs and induced gain functions in Definition \ref{def:induG}. It is then natural to introduce a definition for induced $G$-phases. In Proposition \ref{prop:indu} we will show its consistency with the previous definitions.

\begin{definition}\label{def:induH}
Let $\Gamma =(V_\Gamma, E_\Gamma)$ be a graph and let $H\in \mathcal H_\Gamma$ be a $G$-phase of $\Gamma$.
\begin{itemize}
\item For any $A\subseteq V_\Gamma$,  we define $H_A\in \mathcal H_{\Gamma_A}$ as the  submatrix of $H$ with rows indexed by $A$ and columns indexed by $E_\Gamma(A)$.
\item For any $B\subseteq E_\Gamma$,  we define $H_B\in \mathcal H_{\Gamma_B}$ as the submatrix of $H$ with rows indexed by $V_\Gamma(B)$ and columns indexed by $B$.
\end{itemize}
\end{definition}
When we consider the subgraph $\Gamma_A$ induced by a subset $A$ of $V_\Gamma$, we can assume, without loss of generality,  that the ordering of $V_\Gamma$ and  $E_\Gamma$ are such that $A=\{v_1,\ldots,v_{|A|}\}$ and $E_\Gamma(A)=\{e_1,\ldots,e_{|E_\Gamma(A)|}\}$. In this way, we have $(H_A)_{i,k}=H_{i,k}$ anytime $v_i\in A$ and $e_k\in E_\Gamma(A)$. The analogous assumptions can be made for a subgraph induced by a subset $B$ of $E_\Gamma$. Moreover, we always choose the ordering of the vertices of $L(\Gamma)$ inherited from the ordering of $E_\Gamma$. Thanks to these specifications, we are now in position to prove the next proposition.
\begin{proposition}\label{prop:indu}
Let $H\in \mathcal H_\Gamma$, $\Psi(H)=\psi\in G(\Gamma)$, $\Psi_L(H)=\zeta\in G(L(\Gamma))$. Let $A\subseteq V_\Gamma$ and $B\subseteq E_\Gamma$, so that also $B\subseteq V_{L(\Gamma)}$. Then:
\begin{enumerate}
\item $\Psi(H_A)=\psi_A$;
\item $\Psi(H_B)=\psi_B$;
\item $\Psi_L(H_B)=\zeta_B$.
\end{enumerate}
\end{proposition}
\begin{proof}
Combining Definitions \ref{def:induG}, \ref{def:psi}, \ref{def:induH} we get the following equations.\\
(1) For $v_i,v_j\in A$, $e_k=\{v_i,v_j\}\in E_\Gamma(A)$, using the fact that $\Psi(H)=\psi$, we have:
$$
\psi_A(v_i,v_j)=\psi(v_i,v_j)=s_1 H_{i,k} (H_{j,k})^{-1}=s_1 (H_A)_{i,k} ((H_A)_{j,k})^{-1}=\Psi(H_A)(v_i,v_j).
$$
(2) For $v_i,v_j\in V_\Gamma(B)$, $e_k=\{v_i,v_j\}\in B$, using the fact that $\Psi(H)=\psi$, we have:
$$
\psi_B(v_i,v_j)=\psi(v_i,v_j)=s_1 H_{i,k} (H_{j,k})^{-1}=s_1 (H_B)_{i,k} ((H_B)_{j,k})^{-1}=\Psi(H_B)(v_i,v_j).
$$
(3) For $e_p,e_q\in B$, $v_r=e_p\cap e_q\in V_\Gamma(B)$, using the fact that $\Psi_L(H)=\zeta$, we have:
$$
\zeta_B(e_p,e_q)=\zeta(e_p,e_q)=s_2 (H_{r,p})^{-1} H_{r,q}=s_2 ((H_B)_{r,p})^{-1} (H_B)_{r,q}=\Psi_L(H_B)(e_p,e_q).
$$
\end{proof}

Let $(\Gamma,\psi,H)$ be a triple with $\psi\in G(\Gamma)$, $H\in \mathcal H_\Gamma$ such that $\Psi(H)=\psi$. Such a triple is called \emph{oriented $G$-gain graph} \cite{reff0,line}, and this notion generalizes the one for signed graphs \cite{zaslavskyori,zasmat}. Then the triple $(\Gamma_A,\psi_A,H_A)$ for  $A\subseteq V_\Gamma$ and the triple $(\Gamma_B,\psi_B,H_B)$  for $B\subseteq E_\Gamma$ have the same property: $\Psi(H_A)=\psi_A$ and $\Psi(H_B)=\psi_B$.

Moreover, if $(L(\Gamma),\zeta)$ is a line graph  of the gain graph $(\Gamma, \psi)$, then for every subset $B\subseteq V_{L(\Gamma)}=E_\Gamma$ we have that $(L(\Gamma)_B,\zeta_B)$ is a line graph of the gain subgraph $(\Gamma_B, \psi_B)$.

In light of  Theorem  \ref{theo:vecchio}, the problem of recognizing which gain graphs are gain-line graphs is equivalent to establish the range of the map $\Psi_L$. According to Definition \ref{def:psi}, this range does not depend on the choice of the central involution $s_1$, but only on the choice of $s_2$. For this reason, from now on, we can forget about $s_1$ and we use the notation $s_2=s$.

\begin{remark}\label{rem:K}
From the classical theory we know that the complete graph $K_n$, with $n>3$, is the line graph of the complete bipartite graph $K_{1,n}$, also known as \emph{star graph}, and of no other graph. Notice that a star graph is a particular tree: by Lemma \ref{lem:tree},  whatever the gain function on $K_{1,n}$ is, the associated gain function on its line graph $K_n$ is switching equivalent to $\bold{s}$.
By Proposition \ref{prop:indu}, the same is true for the gain function $\psi_A$ induced by a subset $A$ of the vertices of a gain-line graph $(L,\psi)$ when $L_A$ is isomorphic to $K_n$ for some $n>3$. This provides many examples of gain graphs which are not gain-line graphs, even if the underlying graph is a line graph.
\end{remark}
The previous remark will be crucial for the generalization of Krausz's characterization to gain graphs. However, it is also related to the other characterizations by virtue of the next lemma.

\begin{lemma}\label{lem:completo}
Let $K_n$ be the complete graph on $n$ vertices, with $n\geq 3$. Suppose that $\psi\in G(K_n)$ is such that, for every three distinct vertices $v_0,v_1,v_2\in V_{K_n}$, we have
$$\psi(v_0,v_1)\psi(v_1,v_2)\psi(v_2,v_0)=s.$$ Then $\psi\sim \bold{s}$. In particular,
a gain graph with a complete graph as underlying graph is balanced if and only if all its gain subgraphs on three vertices are balanced.
\end{lemma}
\begin{proof}
By virtue of Lemma \ref{lem:s}, it is enough to prove that, for every closed walk $W$ in $\Gamma$ of length $\ell\geq 3$, one has:
\begin{equation}\label{eq:walk}
\psi(W)= s^\ell=
\begin{cases}
1_G \;&\mbox{ if }  \ell \mbox{ is even}\\
s\; &\mbox{ if }  \ell \mbox{ is odd.}\\
\end{cases}
\end{equation}
We are proving Eq. \eqref{eq:walk} by induction on $\ell$.
Suppose $W$ is a closed walk of length $\ell=3$ visiting vertices: $v_0,v_1,v_2,v_0$.
Clearly $v_0,v_1,v_2$ must be distinct and then Eq. \eqref{eq:walk}
follows from the hypothesis on $\psi$. \\
For a closed walk $W$ of length $\ell$ visiting vertices $v_0,v_1,v_2,v_3,\ldots,v_{\ell-1},v_0$ we consider two cases.
The first is when $v_0=v_2$. In this case, let us define an associated closed walk $W'$ of length $\ell-2$ visiting vertices $v_0,v_3,\ldots v_{\ell-1},v_0$. Clearly $\psi(W)=\psi(W')$ and, for the inductive hypothesis, Eq. \eqref{eq:walk} holds.\\
 In the second case, when $v_0\neq v_2$, consider the associated closed walk $W'$ of length $\ell-1$ visiting vertices $v_0,v_2,\ldots, v_{\ell-1}, v_0$.
 Notice that $v_0,v_1,v_2$ are three distinct vertices and that, by the hypothesis on $\psi$, we have
 $\psi(v_0,v_1)\psi(v_1,v_2)=s\psi(v_0,v_2)$. Therefore:
$$
\psi(W)=\psi(v_0,v_1)\psi(v_1,v_2)\psi(v_2,v_3)\cdots\psi(v_{\ell-1}, v_0)=
 s \psi(W').
$$
Combining with the inductive hypothesis we have proved Eq.  \eqref{eq:walk}.
\end{proof}

Observe that, when $(\Gamma,\psi)$ is a gain graph over an Abelian group $G$, the gain of a closed walk $W$ does not depend on its particular starting vertex $v_0$. This allows to define the gain of an \emph{oriented cycle}. Moreover, the gain of a given oriented cycle does not depend on the particular representative of the switching equivalence class of $\psi$. Actually, it is known that the switching equivalence classes of gain functions on $\Gamma$ are completely determined by their gains on a \emph{cycle basis} (see also the discussion preceding Corollary \ref{cor:cp1}). \\
\indent This argument does not hold in general if the group $G$ is not Abelian: in this case, the choice of two distinct starting vertices produces two conjugate gains. Similarly, two switching equivalent gain functions assign conjugate gains to a given cycle. On the other hand, since $s$ is a central involution of $G$, the property of having gain equal to $s$ or different from $s$ is well defined for a given cycle (regardless its starting vertex and its orientation). For this reason, with a little abuse of notation, for a given triangle subgraph $T$ of $\Gamma$ and a given gain function $\psi\in G(\Gamma)$, we will write $\psi(T)=s$ or $\psi(T)\neq s$. In the same way, even if $G$ is not Abelian, one can check if $\psi \sim {\bf s}$ by only looking at the gains of $\psi$ on a cycle basis. As a consequence, the statement of Lemma \ref{lem:completo} can be reformulated by asking that only the gains of all the triangles of $K_n$ sharing a given vertex $v_0$ are equal to $s$.

\section{Characterizations in gain graphs}

Beineke's characterization of line graphs is in term of  $\mathcal X$-freeness, where  $\mathcal X$ is the list of nine forbidden subgraphs $G_1,\ldots, G_9$  depicted in Fig. \ref{fig:9}. In order to extend it to gain graphs, we introduce the concept of \emph{forbidden induced gain subgraphs}.

\begin{figure}
\begin{center}
\begin{tikzpicture}

\draw (0,0) [black,fill=black] circle (0.08 cm);
\draw (-1,0.8) [black,fill=black] circle (0.08 cm);
\draw (1,0.8) [black,fill=black] circle (0.08 cm);
\draw (0,-0.8) [black,fill=black] circle (0.08 cm);

\draw [-, black] (0,0) -- (-1,0.8);
\draw [-, black] (0,0) -- (1,0.8);
\draw [-, black] (0,0) -- (0,-0.8);
\draw [-, black] (1,0.8) -- (-1,0.8);

\node [align=center] at (-1,-1)
{$P$};

\draw (-0.8+3.5,-0.8) [black,fill=black] circle (0.08 cm);
\draw (+0.8+3.5,-0.8) [black,fill=black] circle (0.08 cm);
\draw (-0.8+3.5,+0.8) [black,fill=black] circle (0.08 cm);
\draw (+0.8+3.5,+0.8) [black,fill=black] circle (0.08 cm);

\draw [-, black] (-0.8+3.5,-0.8) -- (0.8+3.5,-0.8);
\draw [-, black] (-0.8+3.5,-0.8) -- (0.8+3.5,0.8);
\draw [-, black] (-0.8+3.5,-0.8) -- (-0.8+3.5,+0.8);
\draw [-, black] (0.8+3.5,-0.8) -- (0.8+3.5,+0.8);

\draw [-, black] (-0.8+3.5,+0.8) -- (0.8+3.5,+0.8);
\draw [-, black] (-0.8+3.5,+0.8) -- (0.8+3.5,-0.8);

\node [align=center] at (-1+3,-1)
{$K_4$};

\draw [-, black] (0+6.6,0-0.8) -- (1.6+6.6,1.6-0.8);
\draw [-, black] (0+6.6,0-0.8) -- (1.6+6.6,0-0.8);
\draw [-, black] (1.6+6.6,0-0.8) -- (1.6+6.6,1.6-0.8);
\draw [-, black] (0+6.6,0-0.8) -- (0+6.6,1.6-0.8);
\draw [-, black] (0+6.6,1.6-0.8) -- (1.6+6.6,1.6-0.8);

\draw (0+6.6,0-0.8) [black,fill=black] circle (0.08 cm);
\draw (1.6+6.6,0-0.8) [black,fill=black] circle (0.08 cm);
\draw (0+6.6,1.6-0.8) [black,fill=black] circle (0.08 cm);
\draw (1.6+6.6,1.6-0.8) [black,fill=black] circle (0.08 cm);

\node [align=center] at (0+6.6 +0.5,1.6-0.8 -0.5)
{\tiny $T_1$};
\node [align=center] at  (1.6+6.6-0.5,0-0.8+0.5)
{\tiny $T_2$};

\node [align=center] at (-1+7,-1)
{$D$};

\end{tikzpicture}
\end{center}\caption{The Paw graph $P$, the complete graph $K_4$ and the diamond graph $D$.}\label{fig:pkd}
\end{figure}
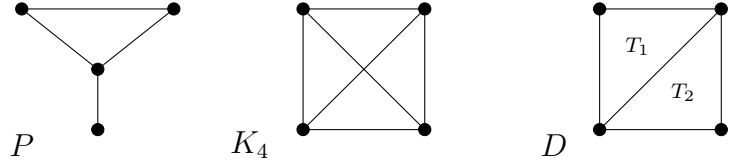

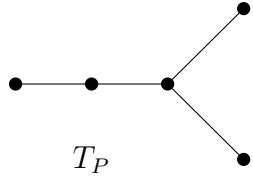
\begin{figure}
\begin{center}
\begin{tikzpicture}

\draw (0,0) [black,fill=black] circle (0.08 cm);
\draw (-1,0) [black,fill=black] circle (0.08 cm);
\draw (-2,0) [black,fill=black] circle (0.08 cm);
\draw (1,1) [black,fill=black] circle (0.08 cm);
\draw (1,-1) [black,fill=black] circle (0.08 cm);

\draw [-, black] (-2,0) -- (-1,0);
\draw [-, black] (-1,0) -- (0,0);
\draw [-, black] (0,0) -- (1,1);
\draw [-, black] (0,0) -- (1,-1);

\node [align=center] at (-1,-1)
{$T_P$};
\end{tikzpicture}
\end{center}\caption{The graph $T_P$ such that $L(T_P)=P$.}\label{fig:t}
\end{figure}

\begin{definition}\label{def:free}
Let $\mathcal Y$ be a set of gain graphs. A gain graph $(\Gamma,\psi)$ is said to be $\mathcal Y$-free if, for any $A\subseteq V_\Gamma$, the gain subgraph $(\Gamma_A,\psi_A)$ is not switching isomorphic to any of the gain graphs in $\mathcal Y$.  The set $\mathcal Y$ is called the list of forbidden gain subgraphs.
\end{definition}
Now we are going to define a set $\mathcal Y$ of forbidden gain subgraphs for the class of gain line graphs.
Not surprisingly, this list include all the gain graphs whose underlying graph is one among $G_1,\ldots, G_9$ of $\mathcal X$ (see Fig. \ref{fig:9}).  We denote it
\begin{equation}\label{eq:xg}
\mathcal X^G:=\{(\Gamma,\zeta) \mid \Gamma \in\mathcal X, \zeta\in G(\Gamma) \}.
\end{equation}
This way, a gain graph $(L,\zeta)$ is $\mathcal X^G$-free if and only if its underlying graph $L$ is $\mathcal X$-free, and so if and only if its underlying graph $L$ is a line graph.
On the other hand, it turns out that this set $\mathcal X^G$ is not big enough in order to characterize, in terms of a list of forbidden gain subgraphs, the class of gain-line graphs. There exist in fact graphs which are line graphs in the classical sense, which become gain-line graphs only when endowed with particular gain functions. For example, it can be easily seen that every gain graph, whose underlying graph is complete, is $\mathcal X^G$-free but it is not necessarily a gain-line graph (see Remark \ref{rem:K}).\\
\indent For this reason, we need to introduce some more gain graphs, whose underlying graphs are  the \emph{paw graph} $P$, the \emph{complete graph} $K_4$ an the \emph{diamond graph} $D$, depicted in Fig. \ref{fig:pkd}. Notice that  the graph $D$ consists of two adjacent triangles, that we denote by $T_1$ and $T_2$. Consider also the graph $T_P$ depicted in Fig. \ref{fig:t}.
We have
\begin{equation}\label{eq:lineclassici}
L(T_P)=P,\qquad L(K_{1,4})=K_4,\qquad L(P)=D.
\end{equation}
Notice that by Whitney's theorem, no other graph has $P$, $K_4$ or $D$ as its line graph. Put:
\begin{equation}\label{eq:y}
\mathcal Y:= \mathcal X^G  \cup  \mathcal F_s
\end{equation}
where $\mathcal X^G$ is defined in Eq. \eqref{eq:xg}, and
\begin{equation}\label{eq:gs}
\mathcal F_s:= \{(P,\zeta) \mid \zeta\nsim \bold{s} \}  \cup  \{(K_4,\zeta) \mid \zeta\nsim \bold{s} \}  \cup \{(D,\zeta) \mid \zeta({T_1})\neq s \mbox{ and } \zeta({T_2})\neq s \}.
\end{equation}
In words, $\mathcal Y$ is the list of all graphs of $\mathcal X$ with every possible gain functions, together with the graphs $P$ and $K_4$ with every gain function non-switching equivalent to $\bold{s}$, together with the diamond graph $D$ with every gain function inducing  gain different from $s$ in both the triangles of $D$.
Notice that every gain graph in $\mathcal Y$ has at most $6$ vertices, exactly as it happened for graphs in $\mathcal X$.
Actually one could consider a narrower list of gain graphs given by a representative of each switching isomorphism class of gain graphs. For example when $G=\mathbb T_2$, the set $\mathcal F_s$ can be replaced by only four signed graphs (see Remark \ref{rem:f}). We can now present the main result of the paper, which is the generalization of Theorem \ref{theo:0} to gain graphs.

\begin{theorem}\label{theo:1}
For a connected gain graph $(L,\zeta)$, with $L=(V_L,E_L)$, the following are equivalent.
\begin{enumerate}
\item $(L,\zeta)$ is a line graph of a gain graph.
\item There exists a partition $E_L=E_1\sqcup E_2 \sqcup\cdots\sqcup E_k$ such that every vertex of $L$ is endpoint of edges from at most two elements of the partition and the induced gain subgraph $(L_{E_i},\zeta_{E_i})$ is a complete graph with $ \zeta_{E_i} \sim \bold{ s}$ in $G(L_{E_i})$, for each $i=1,\ldots,k$.
\item
The following four conditions hold:
\begin{enumerate}[(i)]
\item $L$ is $K_{1,3}$-free;
\item the gain of every odd triangle $T$ of $L$ is $s$;
\item  if $T_1$ and $T_2$ are adjacent odd triangles, then their vertices induce a subgraph  of $L$ that is complete;
\item if $T_1$ and $T_2$ are adjacent even triangles, then the gain of at least one triangle is $s$.
\end{enumerate}
\item $(L,\zeta)$ is  $\mathcal Y$-free.
\item Every gain subgraph of $(L,\zeta)$ induced by a subset of at most $6$ vertices is a line graph of a gain graph.
\end{enumerate}
\end{theorem}
\begin{proof}
${}$\\(2)$\implies$(1)\\
Suppose  that there exists a partition of the edges $E_L=E_1\sqcup E_2 \sqcup\cdots\sqcup E_k$ such that every vertex of $L$ is endpoint of edges from at most two elements of the partition and that the induced gain subgraph $(L_{E_i},\zeta_{E_i})$ is a complete graph with $\zeta_{E_i} \sim \bold{ s}$, for each $i=1,\ldots,k$.
Let us denote by $U=\{u_1,\ldots, u_l\}$ the (possibly empty) subset of vertices of $L$ appearing as endpoints of edges of only one element of the aforementioned partition.\\
\indent Following the classical construction (see for instance  \cite[Theorem~8.4]{Bharary}),
we can define a graph $\Gamma$ whose line graph is $L$. The graph $\Gamma$ is the \emph{intersection graph} of the family of subsets $\{V_L(E_1),\ldots,V_L(E_k),\{u_1\},\ldots,\{ u_l\}\}$ of $V_L$.
Since we need a $G$-phase $H\in \mathcal H_\Gamma$ such that $\Psi_L(H)=\zeta$ (see Theorem \ref{theo:vecchio}), we are going to describe $\Gamma$ with more details.
\\ We set $V_\Gamma=\{x_1,\ldots,x_k, w_1,\ldots,w_l\}$, where the vertices $x_1,\ldots, x_k$ are in a $1$-$1$ correspondence with the parts $E_1,\ldots, E_k$ of the partition of $E_L$, and the vertices $w_1,\ldots, w_l$ are in a $1$-$1$ correspondence with the vertices in $u_1,\ldots, u_l$ of $U$.
We put $x_i\sim x_j$ in $\Gamma$ if $V_L(E_i)\cap V_L(E_j)\neq \emptyset$, and $w_q\sim x_i$ in $\Gamma$
if $u_q\in V_L(E_i)$. Notice that, when $x_i\sim x_j$, the nonempty intersection $V_L(E_i)\cap V_L(E_j)$ consists exactly of one vertex of $L$, because $L_{E_i}$ and $L_{E_j}$ are complete and $E_i$ and $E_j$ are disjoint. Moreover, we set that a vertex $w_i$ is not adjacent to any vertex $w_j$.\\
\indent As stated in \cite[Theorem~8.4]{Bharary}, we have $L(\Gamma)=L$, and  $E_\Gamma =V_L$.
More specifically with an edge $e_p=\{x_i,x_j\} \in E_\Gamma$ is associated the aforementioned unique  vertex of $L$ in $V_L(E_i)\cap V_L(E_j)$  and with an edge  $e_p=\{x_i,w_j\}$ of $\Gamma$ is associated  the vertex $u_j$ of $L$ appearing as endpoint only of edges in $E_i$.
By the hypothesis, for any $i$ we have $\zeta_{E_i}\sim \bold{s}$. By Definition \ref{defswe}, this means that there exists a switching function
$f_i\colon V_L(E_i)\to G$ such that $\zeta_{E_i}^{f_i}=\bold{s}$. We are now able to define the required $H\in \mathcal H_\Gamma$ such that  $\Psi_L(H)=\zeta$.

We first consider the columns of $H$ associated with edges of $\Gamma$ connecting vertices of the type $x_1,\ldots,x_k\in V_\Gamma$.
Suppose we have $e_p=\{x_i,x_j\}\in E_\Gamma$ and then $e_p = V_L(E_i) \cap V_L(E_j)$ as a vertex of $L$.
We set $H_{i,p}:=(f_i(e_p))^{-1}$ and $H_{j,p}:=(f_j(e_p))^{-1}$.
We fill the other entries of the column $p$ with $0$'s.\\
\indent Consider now the columns of $H$ associated with edges of the type $e_p=\{x_i,w_j\}$ of $\Gamma$. In particular $e_p$  is $u_j\in V_L(E_i)$ as a vertex in $L$.
 Then we set $H_{i,p}:=(f_i(e_p))^{-1}$  and $H_{k+j,p}:=1_G$ (or any other group element).	
We fill the other entries of the column $p$ with $0$'s. We can summarize as follows:
\begin{equation}\label{eq:hh}
H_{i,p}:=\begin{cases}
(f_i(e_p))^{-1} &\mbox{if } i\leq k, x_i\in e_p\\
1_G &\mbox{if } i>k, w_{i-k}\in e_p\\
0  &\mbox{otherwise.}
\end{cases}
\end{equation}
By construction, $H\in \mathcal H_\Gamma$. Consider two adjacent vertices $e_p$ and $e_q$ in $L$. They are also incident edges of $\Gamma$,
in particular $e_p\cap e_q$ is not empty. Moreover, $e_p\cap e_q$ must be a vertex in $\{x_1,\ldots,x_k\}\subseteq V_\Gamma$. In fact, if we had
$\{x_i, w_c\},\{x_j,w_c\}\in E_\Gamma$, this would imply that in $L$ the vertex $u_c$ is endpoint of edges of $E_i$ and $E_j$, that is impossible by definition of $U$. Therefore, for any two adjacent vertices $e_p$, $e_q$ of $L$ there must exist $i\in \{1,\ldots,k\}$ such that
 $e_p\cap e_q=x_i$. Combining  Definition \ref{def:psi} and Eq. \eqref{eq:hh}, we have
\begin{equation*}
\begin{split}
\Psi_L(H)(e_p,e_q)&=s (H_{i,p})^{-1} (H_{i,q})= s  (f_i(e_p)^{-1})^{-1} (f_i(e_q))^{-1}= f_i(e_p) s f_i(e_q)^{-1}\\&= \bold{s}^{(f_i)^{-1}}(e_p,e_q)=\zeta_{E_i}(e_p,e_q)=\zeta(e_p,e_q).
\end{split}
\end{equation*}
 It follows that $\zeta=\Psi_L(H)$ and $(L,\zeta)$ is a gain-line graph.

\noindent  (1)$\implies$ (5) \\
Suppose that $(L,\zeta)$ is the gain-line graph of some gain graph $(\Gamma,\psi)$, in particular $V_L=E_\Gamma$.
Then for any subset $A\subseteq V_L$ (in particular, for any subset $A$ with at most $6$ vertices), by virtue of Proposition \ref{prop:indu},  we have that a gain subgraph $(L_A,\zeta_A)$  induced by the subset $A$ of the vertices of $L$ is the gain-line graph of the gain subgraph $(\Gamma_A,\psi_A)$ induced by the subset $A$ of the edges of $\Gamma$.

\noindent (5) $\implies$ (4)\\
If we prove that each gain graph in $\mathcal Y$ (that has at most $6$ vertices)  is not a gain-line graph, then a gain graph $(\Gamma,\psi)$ for which property (5) holds must be
$\mathcal Y$-free (and so property (4) will be satisfied).
Since the underlying graphs of gain graphs in $\mathcal X^G$ are in $\mathcal X$ by Eq. \eqref{eq:xg}, the gain graphs in  $\mathcal X^G$ are not gain-line graphs.
Therefore, it is enough to show that  gain graphs in $\mathcal F_s$ are not gain-line.\\
\indent If $(P,\zeta)$ is the line graph of a gain graph $(\Gamma,\psi)$, the graph $\Gamma$ must be isomorphic to the tree $T_P$ of Fig. \ref{fig:t} (see Eq. \eqref{eq:lineclassici}). By virtue of Lemma \ref{lem:tree} we have  $\zeta\sim \bold{s}$. As a consequence, if  $\zeta\nsim \bold{s}$ then $(P,\zeta)$ cannot be a gain-line graph.\\
Similarly, if $(K_4,\zeta)$ is the gain-line of a graph $(\Gamma,\psi)$, then $\Gamma$ must be isomorphic to the tree $K_{1,4}$ (see Eq. \eqref{eq:lineclassici}) and we can conclude as before.\\
Finally, suppose that  $(D,\zeta)$ is the gain-line of $(\Gamma,\psi)$, so that $\Gamma$ is isomorphic to the paw graph $P$ (see Eq. \eqref{eq:lineclassici}). Consider  $B\subseteq E_\Gamma$ such that $\Gamma_B$ is isomorphic to the graph $K_{1,3}$. We have that $(D_B,\zeta_B)$ is the gain-line of $(\Gamma_B,\psi_B)$, so that $D_B$ is isomorphic to  $K_3$ and, by Lemma \ref{lem:tree}, it must be $\zeta_B\sim \bold{s}$. It follows that, if both the triangles of $(D,\zeta)$ have gains different from $s$, the graph $(D_B,\zeta_B)$ cannot be a subgraph of $(D,\zeta)$ and then  $(D,\zeta)$ cannot be a gain-line graph.

\noindent  (4) $\implies $ (3)\\
Observe that, if $(L,\zeta)$ is $\mathcal Y$-free, then $L$ must be $\mathcal X$-free by definition of $\mathcal{Y}$.
As a consequence of Beineke's characterization, there exists a graph $\Gamma$ such that $L=L(\Gamma)$.
On the other hand, by van Rooij and Wilf's characterization, we have that $L$ is $K_{1,3}$-free and that, if $T_1$ and $T_2$ are two adjacent odd triangles of $L$, then their vertices induce a subgraph of $L$ isomorphic to $K_4$.
Then in order to prove (3) we only have to show that the gain of every odd triangle is $s$  and that when $T_1$ and $T_2$ are adjacent even triangles, at least one among $T_1$ and $T_2$ has gain $s$.\\
We start by proving the first of these two properties.
Suppose that a subset $A:=\{v_1,v_2,v_3\}\subseteq V_L$ induces an odd triangle $T$ in $L$. This implies that there exists $w\in V_L$ such that $|\{v\in A \mid v\sim w \}|$ is $1$ or $3$. Suppose, as a first case, that $|\{v\in A \mid v\sim w \}|=1$. This implies that $w\notin A$ and that
the subgraph of $L$ induced by $A\cup\{ w\}$ is isomorphic to the paw graph $P$. Thanks to conditions (4) we have  $\zeta_{A\cup\{ w\}}\sim \bold{s}$ and then the gain on the triangle is $s$. Suppose now that $|\{v\in A \mid v\sim w \}|=3$. This implies that $w\notin A$ and that
the subgraph of $L$ induced by $A\cup\{ w\}$ is isomorphic to $K_4$. As before, by using condition (4), we have that $\zeta_{A\cup\{ w\}}\sim \bold{s}$ and then the gain of each triangle of this subgraph, and in particular the gain of $T$, is $s$.\\
Let us prove now the second property. Suppose that $T_1$ and $T_2$ are  adjacent even triangles. More precisely, suppose that the vertices of $T_1$ are $\{v_1,v_2,v_3\}$ and the vertices of $T_2$  are  $\{v_1,v_2,v_4\}$. Since $T_1$ and $T_2$ are even, then $v_3\nsim v_4$. In particular, the subgraph of $\Gamma$ induced by $\{v_1,v_2,v_3,v_4\}$ is isomorphic to the diamond graph $D$. By using the condition $(4)$ the gain of at least one of the two triangles must be $s$.

\begin{figure}
\begin{center}
 \begin{tikzpicture}

\fill[mygray](0,0) -- (2,2)--(0,2)--(0,0);

\draw [-, black] (0,0) -- (2,2);
\draw [-, black] (0,0) -- (2,0);
\draw [-, black] (2,0) -- (2,2);
\draw [-, black] (0,0) -- (0,2);
\draw [-, black] (0,2) -- (2,2);

\draw (0,0) [black,fill=black] circle (0.08 cm);
\draw (2,0) [black,fill=black] circle (0.08 cm);
\draw (0,2) [black,fill=black] circle (0.08 cm);
\draw (2,2) [black,fill=black] circle (0.08 cm);

\node [align=center] at (-0.5,-0.5)
{$F_1$};

\fill[mygray](0+4,0) -- (2+4,0)--(5,1)--(0+4,0);
\fill[mygray](0+4,2) -- (2+4,2)--(5,1)--(0+4,2);

\draw [-, black] (0+4,0) -- (2+4,0);
\draw [-, black] (2+4,0) -- (2+4,2);
\draw [-, black] (0+4,0) -- (0+4,2);
\draw [-, black] (0+4,2) -- (2+4,2);

\draw [-, black] (0+4,0) -- (5,1);
\draw [-, black] (2+4,0) -- (5,1);
\draw [-, black] (5,1) -- (0+4,2);
\draw [-, black] (5,1) -- (2+4,2);

\draw (0+4,0) [black,fill=black] circle (0.08 cm);
\draw (2+4,0) [black,fill=black] circle (0.08 cm);
\draw (0+4,2) [black,fill=black] circle (0.08 cm);
\draw (2+4,2) [black,fill=black] circle (0.08 cm);
\draw (1+4,1) [black,fill=black] circle (0.08 cm);

\node [align=center] at (-0.5+4,-0.5)
{$F_2$};
  \fill[mygray]( 1.2+8,1)-- (0+8,0)--(1.5+8,2)--( 1.2+8,1);
 \fill[mygray]( 1.5+8,0.4)-- (0+8,0)--(3+8,0)--( 1.5+8,0.4);
\fill[mygray]( 1.8+8,1)--( 1.2+8,1)--( 1.5+8,0.4)--( 1.8+8,1);
   \fill[mygray]( 1.8+8,1)--(1.5+8,2)--(3+8,0)--( 1.8+8,1);

 \draw [-, black] (0+8,0) -- (3+8,0);
\draw [-, black] (0+8,0) -- (1.5+8,2);
\draw [-, black] (1.5+8,2) -- (3+8,0);

 \draw [-, black]( 1.5+8,0.4)-- (0+8,0);
 \draw [-, black]( 1.5+8,0.4)-- (3+8,0);

\draw [-, black]( 1.2+8,1)-- (0+8,0);
\draw [-, black]( 1.2+8,1)-- (1.5+8,2);

\draw [-, black]( 1.8+8,1)--(1.5+8,2);
 \draw [-, black]( 1.8+8,1)--(3+8,0);

\draw [-, black]( 1.8+8,1)--( 1.2+8,1);
   \draw [-, black]( 1.8+8,1)--( 1.5+8,0.4);
     \draw [-, black]( 1.5+8,0.4)--( 1.2+8,1);

\draw (0+8,0) [black,fill=black] circle (0.08 cm);
\draw (3+8,0) [black,fill=black] circle (0.08 cm);
\draw (1.5+8,2) [black,fill=black] circle (0.08 cm);

\draw( 1.8+8,1)  [black,fill=black] circle (0.08 cm);
\draw( 1.2+8,1)  [black,fill=black] circle (0.08 cm);
\draw( 1.5+8,0.4)  [black,fill=black] circle (0.08 cm);

\node [align=center] at (-0.5+8,-0.5)
{$F_3$};
\end{tikzpicture}
\end{center}\caption{The exceptional graphs $F_1, F_2, F_3$.}\label{fig:tre}
\end{figure}
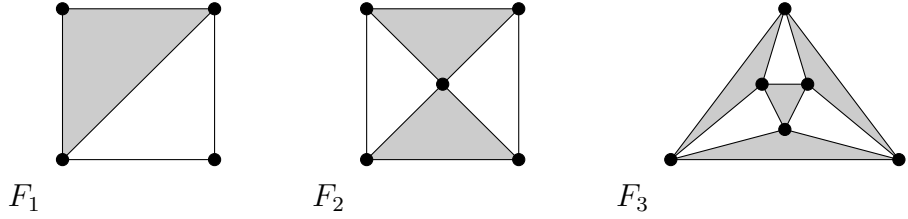

\noindent(3) $\implies$ (2)\\
As in the classical setting, we consider two distinct cases, when $L$ contains   adjacent  even triangles and when it does not.\\
In \cite{beineke} it is shown  that there exist exactly three exceptional graphs containing two adjacent even triangles and satisfying conditions (i) and (iii) of (3): they are the graphs $F_1$, $F_2$, $F_3$, depicted in Fig. \ref{fig:tre}. The last condition on $(L,\zeta)$ in (3) about adjacent even triangles ensures that, if $L$ is either $F_1$, $F_2$, or $F_3$, we can assume that  at least the gains  of all triangles filled in gray are $s$, or  at least the gains  of all triangles unfilled, are $s$.
Let us assume to be in the first case, so that the  gain of each gray triangle is $s$ (the second case can be
similarly covered). We define a partition of $E_L$ in the following way:
the edges bounding the same gray triangle are in the same part; each of the other parts contains exactly one of the (possibly) remaining edges.
This way, the  subgraph induced by each part is isomorphic to $K_2$ or $K_3$. Moreover, the induced gain function on this subgraph, in both cases,  is switching equivalent to $\bold{s}$. In fact, if the subgraph is isomorphic to $K_3$, this is true because we assumed that the gain of each gray triangle is $s$; if it is isomorphic to $K_2$, this is true because  $K_2$ is a tree. Finally, from the picture, it is clear that a vertex is endpoint of edges from at most two parts. \\
\indent Suppose now that in $L$ there is no pair of adjacent even triangles. Then, according to the construction in the proof of the main theorem in \cite{beineke}, there exists a partition $E_L=E_1\sqcup E_2\sqcup \cdots \sqcup E_k$ such that:
\begin{itemize}
\item $\Gamma_{E_i}$ is complete;
\item if $|E_i|=3$, then the endpoints of $E_i$ form an odd triangle;
\item every vertex $v\in V_L$ is endpoint of edges from at most two parts.
\end{itemize}
Then if $E_i$ consists of exactly one edge, clearly $\zeta_{E_i}\sim \bold{s}$. If $|E_i|=3$, the associated triangle is odd and, by hypothesis (3), we have $\zeta_{E_i}\sim \bold{s}$. Finally suppose that $\Gamma_{E_i}$  is isomorphic to $K_n$ with $n>3$. Notice that every triangle of $K_n$ is odd and then, the gain of each triangle is $s$. By Lemma \ref{lem:completo} we have $\zeta_{E_i}\sim \bold{s}$ and we have done.
\end{proof}

\begin{example}\label{examplenuovo}
Consider the $G$-gain graph $(L,\xi)$ depicted in Fig.~\ref{fig:nuova1}, where $G$ is any group  and $a,b,c,d\in G$.

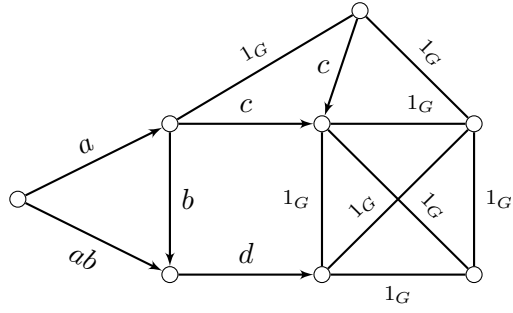
\begin{figure}
\begin{center}

\begin{tikzpicture}
% vertices
\node[vertex]  (1) at (0,0){};
\node[vertex]  (2) at (2,1) {};
\node[vertex]  (3) at (2,-1) {};

\node[vertex]  (4) at (4,-1) {};
\node[vertex]  (5) at (4,1) {};
\node[vertex]  (6) at (6,-1) {};
\node[vertex]  (7) at (6,1) {};
\node[vertex]  (8) at (4.5,2.5) {};

%edges
\draw[dedge] (1) -- (2) node[midway, above,sloped] {\small $a$};
\draw[dedge] (2) -- (3) node[midway, right] {\small $b$};
\draw[dedge] (1) -- (3)node[midway, below,sloped] {\small $ab$};

\draw[dedge] (2) -- (5) node[midway, above,sloped] {\small $c$};
\draw[dedge] (8) -- (5) node[midway,left] {\small $c$};
\draw[edge] (2) -- (8)node[midway, above,sloped] {\tiny $1_G$};

\draw[edge] (8) -- (7)node[midway, above,sloped] {\tiny $1_G$};

\draw[dedge] (3) -- (4) node[midway, above,sloped] {\small $d$};

\draw[edge] (4) -- (5) node[midway, left] {\tiny $1_G$};
\draw[edge] (4) -- (6)node[midway, below,sloped] {\tiny $1_G$};
\draw[edge] (4) -- (7)node[pos=0.35, above,sloped] {\tiny $1_G$};
\draw[edge] (5) -- (6)node[pos=0.65, above,sloped] {\tiny $1_G$};
\draw[edge] (5) -- (7)node[pos=0.67, above,sloped] {\tiny $1_G$};
\draw[edge] (6) -- (7)node[midway, right] {\tiny $1_G$};
\end{tikzpicture}
\end{center}\caption{The $G$-gain graph $(L,\xi)$ of Example \ref{examplenuovo}.}\label{fig:nuova1}
\end{figure}
We will omit the label $1_G$ on each edge whose gain is $1_G$. At first we observe in Fig.~\ref{fig:nuova2} that the underlying graph $L$ admits the following partition of the edge set
$$E_L=E_1\sqcup E_2 \sqcup E_3 \sqcup E_4 \sqcup E_5,$$
such that every vertex in $V_L$ is endpoint of edges from at most two elements of the partition and the subgraph $L_{E_i}$ of $L$ induced by $E_i$ is a complete graph, for each $i=1,\ldots, 5$.
By Krausz's characterization in Theorem \ref{theo:0}, the graph $L$ is the line graph of some simple graph. Moreover, one can check that the edge-induced gain subgraph $(L_{E_i},\zeta_{E_i})$ is balanced, that is $ \zeta_{E_i} \sim \bold{ 1_G}$ in $G(L_{E_i})$, for each $i=1,\ldots,5$.
Hence the condition $(2)$ of Theorem \ref{theo:1} holds for $s=1_G$, and then  $(L,\xi)$ is a gain-line graph. We are going to explicitly construct
the graph $\Gamma$ such that $L(\Gamma)=L$ and a $G$-phase $H\in \mathcal H_\Gamma$ such that $\psi_L(H)=\xi$ (see Definition \ref{def:psi}).

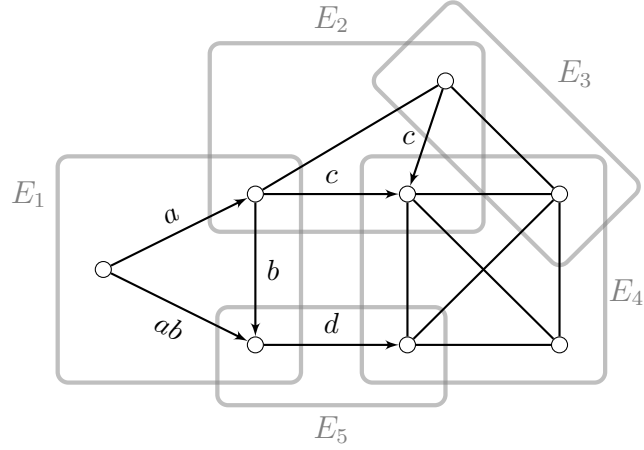
\begin{figure}
\begin{center}
\begin{tikzpicture}
\draw[ultra thick, gray,opacity=0.5, rounded corners] (-0.6,-1.5) rectangle (2.6,1.5);
\node[gray,left] at (-0.6,1){$E_1$};

\draw[ultra thick, gray,opacity=0.5,rounded corners] (1.4,0.5) rectangle (5,3);
\node[gray,above] at (3,3){$E_2$};

\draw[ultra thick, gray,opacity=0.5,rounded corners] (1.5,-0.5) rectangle (4.5,-1.8);
\node[gray,below] at (3,-1.8){$E_5$};

\draw[ultra thick, gray,opacity=0.5,rounded corners] (3.4,-1.5) rectangle (6.6,1.5);
\node[gray,right] at (6.5,-0.3){$E_4$};

\draw[ultra thick, gray,opacity=0.5,rounded corners] (3.5,2.5)--(4.6,3.6)--(7.1,1.1)--(6,0)--cycle;
\node[gray] at (6.2,2.6){$E_3$};
% vertices
\node[vertex]  (1) at (0,0){};
\node[vertex]  (2) at (2,1) {};
\node[vertex]  (3) at (2,-1) {};

\node[vertex]  (4) at (4,-1) {};
\node[vertex]  (5) at (4,1) {};
\node[vertex]  (6) at (6,-1) {};
\node[vertex]  (7) at (6,1) {};
\node[vertex]  (8) at (4.5,2.5) {};

%edges
\draw[dedge] (1) -- (2) node[midway, above,sloped] {\small $a$};
\draw[dedge] (2) -- (3) node[midway, right] {\small $b$};
\draw[dedge] (1) -- (3)node[midway, below,sloped] {\small $ab$};

\draw[dedge] (2) -- (5) node[midway, above,sloped] {\small $c$};
\draw[dedge] (8) -- (5) node[midway,left] {\small $c$};
\draw[edge] (2) -- (8)node[midway, above,sloped] {\tiny };

\draw[edge] (8) -- (7)node[midway, above,sloped] {\tiny };

\draw[dedge] (3) -- (4) node[midway, above,sloped] {\small $d$};

\draw[edge] (4) -- (5) node[midway, left] {\tiny };
\draw[edge] (4) -- (6)node[midway, below,sloped] {\tiny };
\draw[edge] (4) -- (7)node[pos=0.35, above,sloped] {\tiny };
\draw[edge] (5) -- (6)node[pos=0.65, above,sloped] {\tiny };
\draw[edge] (5) -- (7)node[pos=0.67, above,sloped] {\tiny };
\draw[edge] (6) -- (7)node[midway, right] {\tiny };

\end{tikzpicture}
\end{center}\caption{The partition of the edge-set  of $(L,\xi)$ of Example \ref{examplenuovo}.}\label{fig:nuova2}
\end{figure}
Following the formalism of the proof of Theorem \ref{theo:1} (implication (2)$\implies$(1)), we denote by $u_1,u_2$ the two vertices of $L$ that are endpoints of edges of only one element of the aforementioned partition, and by $v_1,\ldots,v_6$ the other vertices. The graph $\Gamma =(V_\Gamma, E_\Gamma)$ can be described as the intersection graph of the family of subsets $\{V_L(E_1),\ldots,V_L(E_5),\{u_1\},\{ u_2\}\}$. More precisely, we set $V_\Gamma=\{x_1,\ldots,x_5,w_1,w_2\}$ in such a way that
\begin{enumerate}
\item $x_1$ is associated  with $V_L(E_1)=\{u_1,v_1,v_2\}$;
\item $x_2$ is associated  with $V_L(E_2)=\{v_1,v_4,v_5\}$;
\item $x_3$ is associated  with $V_L(E_3)=\{v_5,v_6\}$;
\item $x_4$ is associated  with $V_L(E_4)=\{v_3,v_4,v_6,u_2\}$;
\item $x_5$ is associated  with $V_L(E_5)=\{v_2,v_3\}$;
\item $w_1$ is associated  with the singleton $\{u_1\}$;
\item $w_2$ is associated  with the singleton $\{u_2\}$.
\end{enumerate}
Moreover, each edge of $\Gamma$ is associated with the only vertex of $L$ that is in the intersection of its endpoints. See Fig.~\ref{fig:nuova3}.
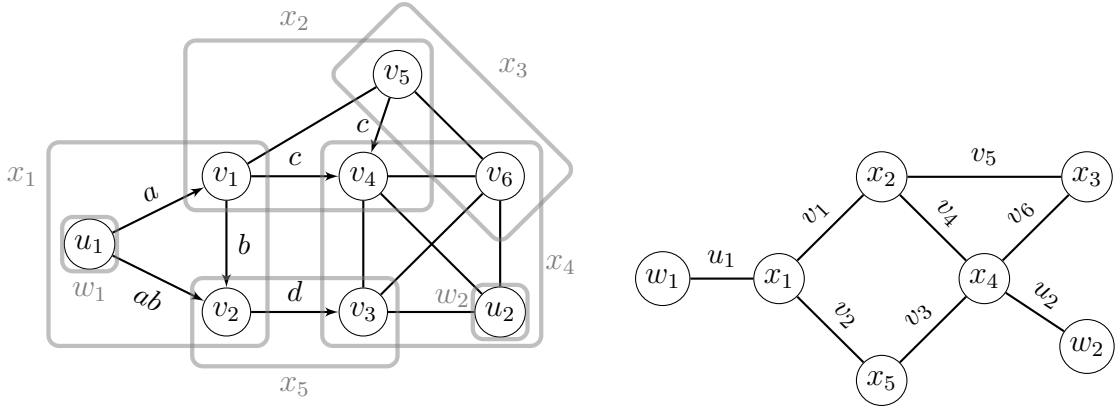
\begin{figure}[h]
\begin{center}
\begin{tikzpicture}[scale=0.90]
% vertices
\node[vertex]  (1) at (0,0){$u_1$};
\node[vertex]  (2) at (2,1) {$v_1$};
\node[vertex]  (3) at (2,-1) {$v_2$};

\node[vertex]  (4) at (4,-1) {$v_3$};
\node[vertex]  (5) at (4,1) {$v_4$};
\node[vertex]  (6) at (6,-1) {$u_2$};
\node[vertex]  (7) at (6,1) {$v_6$};
\node[vertex]  (8) at (4.5,2.5) {$v_5$};

%edges
\draw[dedge] (1) -- (2) node[midway, above,sloped] {\small $a$};
\draw[dedge] (2) -- (3) node[midway, right] {\small $b$};
\draw[dedge] (1) -- (3)node[midway, below,sloped] {\small $ab$};

\draw[dedge] (2) -- (5) node[midway, above,sloped] {\small $c$};
\draw[dedge] (8) -- (5) node[midway,left] {\small $c$};
\draw[edge] (2) -- (8)node[midway, above,sloped] {\tiny };

\draw[edge] (8) -- (7)node[midway, above,sloped] {\tiny };

\draw[dedge] (3) -- (4) node[midway, above,sloped] {\small $d$};

\draw[edge] (4) -- (5) node[midway, left] {\tiny };
\draw[edge] (4) -- (6)node[midway, below,sloped] {\tiny };
\draw[edge] (4) -- (7)node[pos=0.35, above,sloped] {\tiny };
\draw[edge] (5) -- (6)node[pos=0.65, above,sloped] {\tiny };
\draw[edge] (5) -- (7)node[pos=0.67, above,sloped] {\tiny };
\draw[edge] (6) -- (7)node[midway, right] {\tiny };

\draw[ultra thick, gray,opacity=0.5, rounded corners] (-0.6,-1.5) rectangle (2.6,1.5);
\node[gray,left] at (-0.6,1){$x_1$};

\draw[ultra thick, gray,opacity=0.5,rounded corners] (1.4,0.5) rectangle (5,3);
\node[gray,above] at (3,3){$x_2$};

\draw[ultra thick, gray,opacity=0.5,rounded corners] (1.5,-0.5) rectangle (4.5,-1.8);
\node[gray,below] at (3,-1.8){$x_5$};

\draw[ultra thick, gray,opacity=0.5,rounded corners] (3.4,-1.5) rectangle (6.6,1.5);
\node[gray,right] at (6.5,-0.3){$x_4$};

\draw[ultra thick, gray,opacity=0.5,rounded corners] (3.5,2.5)--(4.6,3.6)--(7.1,1.1)--(6,0)--cycle;
\node[gray] at (6.2,2.6){$x_3$};

\draw[ultra thick, gray,opacity=0.5,rounded corners] (-0.4,-0.4) rectangle (0.4,0.4);
\node[gray,below] at (0,-0.4){$w_1$};

\draw[ultra thick, gray,opacity=0.5,rounded corners] (-0.4+6,-0.4-1) rectangle (0.4+6,0.4-1);
\node[gray,below] at (5.3,-0.5){$w_2$};

\end{tikzpicture}\hspace{5mm}
\begin{tikzpicture}[scale=0.90]

% vertices
\node[vertex]  (x1) at (0,0){$x_1$};
\node[vertex]  (x2) at (1.5,1.5) {$x_2$};
\node[vertex]  (x4) at (3,0) {$x_4$};
\node[vertex]  (x5) at (1.5,-1.5) {$x_5$};
\node[vertex]  (w1) at (-1.7,0){$w_1$};
\node[vertex]  (w2) at (4.5,-1){$w_2$};
\node[vertex]  (x3) at (4.5,1.5){$x_3$};

%edges

\draw[edge] (x1) -- (x2) node[midway, above,sloped] {\small $v_1$};
\draw[edge] (x2) -- (x4)node[midway, above,sloped] {\small $v_4$};
\draw[edge] (x4) -- (x5)node[midway, above,sloped] {\small $v_3$};
\draw[edge] (x5) -- (x1)node[midway, above,sloped] {\small $v_2$};
\draw[edge] (w1) -- (x1) node[midway, above,sloped] {\small $u_1$};
\draw[edge] (w2) -- (x4) node[midway, above,sloped] {\small $u_2$};
\draw[edge] (x3) -- (x4) node[midway, above,sloped] {\small $v_6$};
\draw[edge] (x3) -- (x2) node[midway, above,sloped] {\small $v_5$};

\end{tikzpicture}

\end{center}\caption{The family of subsets of vertices of $(L,\xi)$ and its intersection graph isomorphic to  $\Gamma$  of Example \ref{examplenuovo}.}\label{fig:nuova3}
\end{figure}

\noindent For each $i=1,\ldots, 5$, let us define a switching function $f_i$ on the vertex set $V_L(E_i)$ as follows:
\begin{itemize}
\item $f_1\colon V_L(E_1)\to G$ such that
$$
f_1(v_1)=1_G,\quad f_1(v_2)=b^{-1},\quad f_1(u_1)=a;
$$
\item $f_2\colon V_L(E_2)\to G$  such that
$$
f_2(v_1)=1_G,\quad f_2(v_4)=c^{-1},\quad f_2(v_5)=1_G;
$$
\item $f_3\colon V_L(E_3)\to G$  such that
$$
f_3(v_5)=f_3(v_6)=1_G;
$$
\item $f_4\colon V_L(E_4)\to G$  such that
$$
f_4(v_3)=f_4(v_4)=f_4(v_6)=f_4(u_2)=1_G;
$$
\item $f_5\colon V_L(E_5)\to G$  such that
$$
f_5(v_2)=d,\quad f_5(v_3)=1_G.
$$
\end{itemize}
One can directly check that the condition $\xi_{E_i}^{f_i}=\bold{1_G}$ is satisfied for each $i$. We can summarize the values of these functions in this table:
\begin{center}
\begin{tabular}{c|c|c|c|c|c|c|c|c|}
& $v_1$ & $v_2$ & $v_3$ & $v_4$ & $v_5$ & $v_6$ & $u_1$ & $u_2$\\
\hline
$f_1$ & $1_G$  & $b^{-1}$  & $\diagup$ & $\diagup$ &$\diagup$& $\diagup$&$a$ &$\diagup$ \\
\hline
$f_2$ & $1_G$& $\diagup$  & $\diagup$ & $c^{-1}$ &$1_G$& $\diagup$&$\diagup$&$\diagup$ \\
\hline
$f_3$ & $\diagup$ & $\diagup$ & $\diagup$ & $\diagup$  &$1_G$ & $1_G$&$\diagup$ &$\diagup$ \\
\hline
$f_4$ & $\diagup$ & $\diagup$ & $1_G$ & $1_G$ &$\diagup$ & $1_G$&$\diagup$ &$1_G$\\
\hline
$f_5$ & $\diagup$ & $d$ & $1_G$  & $\diagup$  &$\diagup$ & $\diagup$&$\diagup$ &$\diagup$\\
\hline
\end{tabular}
\end{center}
Finally, according with Eq. \eqref{eq:hh}, we obtain the $G$-phase
$$
H=\begin{pmatrix}
1_G & b & 0 &0 & 0& 0 &a^{-1}&0\\
1_G & 0 & 0 &c & 1_G& 0 &0 &0\\
0 & 0 & 0 &0 & 1_G& 1_G &0 &0\\
0 & 0 & 1_G & 1_G &0 & 1_G& 0 &1_G \\
0 & d^{-1} &1_G & 0 &0 & 0& 0 &0 \\
0 & 0 &0 & 0 &0 & 0& 1_G& 0 \\
0 & 0 &0 & 0 &0 & 0& 0& 1_G \\
\end{pmatrix}.
$$
Notice that the first five rows can be obtained from the previous table by taking the inverse of the entries of the table, and by replacing the empty spaces with $0\in\mathbb CG$. One can check that $\psi_L(H)=\xi$.
For instance, one has:
$$
\psi_L(H)(u_1,v_2)= H_{1,7}^ {-1}H_{1,2}=ab = \xi(u_1,v_2); \quad
\psi_L(H)(v_2,v_3)= H_{5,2}^ {-1}H_{5,3}=d=\xi(v_2,v_3).
$$
\end{example}
\vspace{1cm}

Under the hypothesis that $L=L(\Gamma)$, that is, $L$ is a line graph as underlying graph, it is very easy to establish whether a gain graph $(L,\zeta)$ is a gain-line or not, as the next corollary shows.
\begin{corollary}\label{cor:primo}
Suppose $(L,\zeta)$ is such that $L=L(\Gamma)$ and $L$ is not isomorphic to either $F_1$, $F_2$, or $F_3$ of Fig. \ref{fig:tre}.
Then $(L,\zeta)$ is a gain-line if and only if each of its odd triangles has gain $s$.
\end{corollary}
\begin{proof}
We are using the characterization (3) of Theorem \ref{theo:1}. Clearly if $(L,\zeta)$ is a gain-line, then every odd triangle has gain $s$. Let us prove the converse implication. Since $L=L(\Gamma)$, by the van Rooij and Wilf's characterization, we know that $L$ is $K_{1,3}$-free (condition (i) of (3)) and that if $T_1$ and $T_2$ are adjacent odd triangles, then their vertices induce a subgraph  of $L$ that is complete (condition (iii) of (3)). Moreover, since $L$ is not isomorphic to either $F_1$, $F_2$, or $F_3$, it has no adjacent even triangles \cite{beineke}, and the condition (iv) of (3) is trivially satisfied. Since the gain of every odd triangle $T$ of $L$ is $s$ by the hypothesis (condition (ii) of (3)), the proof is completed.
\end{proof}
As a consequence of Corollary \ref{cor:primo} we have that, if $(L,\zeta)$ is a gain graph such that:
\begin{itemize}
\item either there is an odd triangle $T$ with gain that is not central,
\item or there is an odd triangle $T$ with gain that is not an involution,
\item or there are odd triangles $T_1,T_2$ with different gains,
\end{itemize}
then $(L,\zeta)$ cannot be a gain-line, for any choice of $s$.

\begin{remark}\label{rem25nov}\rm
When $s=1_G$,  if $(L,\zeta)$ is a gain-line graph, then all its odd triangles are balanced. Under the hypothesis that $L=L(\Gamma)$ and that $L$ is not isomorphic to either $F_1$, $F_2$ or $F_3$, the condition of balance of all odd triangles of $(L,\zeta)$  is also sufficient to say that $(L,\zeta)$ is a gain-line.\\
\indent When $G=\mathbb T$ and $s=-1$, if $L$  is a line graph which is not isomorphic to either $F_1$, $F_2$ or $F_3$, then we have that a $\mathbb T$-gain graph $(L,\zeta)$ is a gain-line
 if and only if all odd triangles are \emph{antibalanced} (i.e., their gains are equal to $-1$).
\end{remark}
Starting from a line graph $L(\Gamma)$, and using a $G$-phase $H\in \mathcal{H}_\Gamma$, it is always possible to define a gain $\zeta\in G(L(\Gamma))$ such that $(L(\Gamma),\zeta)$ is a gain graph (see Remark \ref{rem:s2}). On the other hand, in many cases there exist gain functions $\psi$ such that $(L(\Gamma),\psi)$ is not a gain-line graph (e.g., this happens for the graphs from the family $\mathcal{F}_s$). In this setting, there are two families of graphs with a special behavior, as the next corollary shows.
\begin{corollary}\label{cor:cp}
Let $G$ be a nontrivial group. Then a connected graph $L=(V_L,E_L)$  is such that  $(L,\zeta)$ is a gain-line for all  $\zeta\in G(L)$ (that is, the map $\mathcal L$ is surjective) if and only if $L$ is a path or a cycle.
\end{corollary}
\begin{proof}
If $L$ is a path or a cycle the condition of Corollary \ref{cor:primo}  is trivially satisfied for any $\zeta\in G(L)$. \\
\indent Let us prove the converse implication. First of all, notice that for $|V_L|\leq 3$ the claim is trivial. Therefore, we can assume $|V_L|\geq 4$. By contradiction, suppose that $(L,\zeta)$ is a gain-line for any $\zeta\in G(L)$ but $L$ is not a path nor a cycle.
As a first step we prove that, under these hypotheses, $L$ contains at least one odd triangle or at least a pair of adjacent triangles.
\\ The hypothesis that  $L$ is not a path nor a cycle implies that  there exists a vertex $v_0\in V_L$ that is adjacent to three distinct vertices $v_1,v_2,v_3$.
Since $(L,\zeta)$ is a gain-line  for every $\zeta$, then in particular the underlying graph $L$ is a line graph and it is $K_{1,3}$-free (see Theorem \ref{theo:0}). As a consequence,  at least two among $v_1,v_2,v_3$  are adjacent. Suppose, without loss of generality, that $v_1\sim v_2$. Now,
if $v_3$ is adjacent to neither $v_1$ nor to $v_2$, the triangle $v_0,v_1,v_2$ is odd, but if $v_3$ is adjacent to at least one among $v_1$ and $v_2$, then there are two adjacent triangles.\\
In the first case, $L$ has an odd triangle $T$. Since $G$ is nontrivial, it is possible to define a gain function $\zeta_1$ such that $\zeta_1(T)\neq s$. By Corollary \ref{cor:primo}, the gain graph $(L,\zeta_1)$ is not a gain-line, which is a contradiction.
In the second case, if at least one of the  two adjacent triangles is odd, we can argue as in the previous case. However, if both the adjacent triangles $T_1$ and $T_2$ are even, it is always possible to define a gain function $\zeta_2$ such that $\zeta_2(T_1)\neq s$ and $\zeta_2(T_2)\neq s$. By Theorem \ref{theo:1}, the gain graph
$(L,\zeta_2)$ is not a gain-line, a contradiction again.
\end{proof}

\section{Spectral characterizations in signed graphs}

In this section we assume that $G=\mathbb T_2=\{\pm 1\}$. Then a $\mathbb T_2$-gain graph $(\Gamma,\sigma)$ is usually called a signed graph, the gain function $\sigma$ is usually known as the signature of $\Gamma$, the gain of a walk is usually called the sign of a walk.
Also the cycles can be partitioned into positive (balanced) cycles and negative (unbalanced) cycles, according to the parity of the number of negative edges.
By fixing an ordering $\{v_1,\ldots,v_n\}$ of the vertices of $\Gamma$, it is possible to define an \emph{adjacency matrix of the signed graph $(\Gamma,\sigma)$}, denoted with  $A_{(\Gamma,\sigma)}$, by setting:
$$
\left(A_{(\Gamma,\sigma)}\right)_{i,j}:=\begin{cases}
\sigma(v_i,v_j) &\mbox{ if } v_i\sim v_j\\
0 &\mbox{otherwise.}
\end{cases}
$$
By definition, the matrix $A_{(\Gamma,\sigma)}$ is real and symmetric. The \emph{spectrum $\Spec(\Gamma,\sigma)$ of the signed graph} is defined as the spectrum of the matrix $A_{(\Gamma,\sigma)}$ and it is invariant under switching isomorphism. Notice that, with any signature $\sigma$ of $\Gamma$, the opposite signature $-\sigma$ is associated, and it satisfies $\Spec(\Gamma,-\sigma)=-\Spec(\Gamma,\sigma)$.

For a connected graph $\Gamma$ with $n$ vertices and $m$ edges there exists a spanning tree $T$ of $\Gamma$ that is a subgraph induced by $n-1$ of its edges. As we have already seen in the proof of Lemma \ref{lem:s}, each of the remaining $m-n+1$ edges of $\Gamma$ induces exactly one cycle involving only edges from $T$ and the edge itself. This set of $m-n+1$ cycles is  a \emph{cycle basis}
(see, for example, \cite{Bharary}), and the number $m-n+1$ is called \emph{circuit rank} of $\Gamma$. A standard modification of the argument in the proof of Lemma \ref{lem:s} proves that two signatures are switching equivalent if and only if their
 signs on a cycle basis coincide. In particular, the cardinality of the  switching equivalence  classes of signatures $[\mathbb T_2(\Gamma)]$ on a connected graph $\Gamma$ is $2^{m-n+1}$ (see \cite[Proposition~3.1]{nase}). Since $\mathcal L$ is an injective map from $[\mathbb T_2(\Gamma)]$ to $[\mathbb T_2(L(\Gamma))]$, the characterization of surjectivity of $\mathcal L$  in Corollary \ref{cor:cp} gives information on the circuit rank of the line graph.

\begin{corollary}\label{cor:cp1}
If $\Gamma$ is neither a path nor a cycle, the circuit rank of $L(\Gamma)$ is greater than the circuit rank of $\Gamma$.
\end{corollary}
Now we are going to analyze the line construction in signed graphs in relation with the two possible choices of $s$.\\
\indent Suppose $s=1$. As observed in Remark \ref{rem:s2},
the signed graph $(\Gamma,\bold{+1})$ is a line graph if and only if $\Gamma$ is a line graph.
 Notice that the usual way to embed an unsigned graph $\Gamma$ into signed graphs, is through the all-positive signature $\bold{+1}$:
 $$\Gamma\mapsto (\Gamma,\bold{+1}).$$
As a consequence, a graph that is a line graph, regarded as a signed graph with the all-positive signature, is still a line graph.\\
\indent On the other hand, if $s=-1$, we have  that the signed graph $(\Gamma,\bold{-1})$ is a line graph if and only if $\Gamma$ is a line graph.
As we have shown in Corollary \ref{cor:primo}, if $(L,\sigma)$ is such that $L=L(\Gamma)$ and $L$ is not isomorphic to either $F_1$, $F_2$, or $F_3$, the signed graph $(L,\sigma)$ is a line graph with $s=1$ if and only if all its odd triangles are balanced; $(L,\sigma)$ is a line graph with $s=-1$ if and only all its odd triangles are unbalanced. In particular, if a signed graph has simultaneously balanced and unbalanced odd triangles, it cannot be a signed line graph, with any choice of $s$.

In what follows we want to give a spectral characterization of the class of line graphs of signed graphs. For this reason, we look again at the set of the forbidden signed subgraphs.
\begin{remark}\label{rem:f}
It is not difficult to check that every signed graph in $\mathcal F_1$ (defined in Eq. \eqref{eq:gs}) is switching isomorphic to one among $(P,\bold{-1})$, $(K_4,\bold{-1})$, $(K_4,\sigma_1)$, $(D,\bold{-1})$,  where  $(K_4,\sigma_1)$ is a complete signed graph  with exactly one negative edge (see Fig. \ref{fig:fsegn}).
 Similarly, every signed graph in $\mathcal F_{-1}$  is switching isomorphic to one among $(P,\bold{1})$, $(K_4,\bold{1})$, $(K_4,\sigma_1)$, $(D,\bold{1})$ (see Fig. \ref{fig:fsegnmeno}).
 \end{remark}

\begin{figure}[h]
\begin{center}
 \begin{tikzpicture}

\draw (-0.1+0,0) [black,fill=black] circle (0.08 cm);
\draw (-0.1-1,0.8) [black,fill=black] circle (0.08 cm);
\draw (-0.1+1,0.8) [black,fill=black] circle (0.08 cm);
\draw (-0.1+0,-0.8) [black,fill=black] circle (0.08 cm);

\draw [-,dashed,  black] (-0.1+0,0) -- (-0.1-1,0.8);
\draw [-,dashed,  black] (-0.1+0,0) -- (-0.1+1,0.8);
\draw [-,dashed,  black] (-0.1+0,0) -- (-0.1+0,-0.8);
\draw [-,dashed, black] (-0.1+1,0.8) -- (-0.1-1,0.8);

\node [align=center] at (-0.1-1+0.15,-1)
{$(P,\bold{-1})$};

\draw (-0.8+3.5,-0.8) [black,fill=black] circle (0.08 cm);
\draw (+0.8+3.5,-0.8) [black,fill=black] circle (0.08 cm);
\draw (-0.8+3.5,+0.8) [black,fill=black] circle (0.08 cm);
\draw (+0.8+3.5,+0.8) [black,fill=black] circle (0.08 cm);

\draw [-,dashed,  black] (-0.8+3.5,-0.8) -- (0.8+3.5,-0.8);
\draw [-, dashed, black] (-0.8+3.5,-0.8) -- (0.8+3.5,0.8);
\draw [-, dashed, black] (-0.8+3.5,-0.8) -- (-0.8+3.5,+0.8);
\draw [-, dashed, black] (0.8+3.5,-0.8) -- (0.8+3.5,+0.8);

\draw [-,dashed,  black] (-0.8+3.5,+0.8) -- (0.8+3.5,+0.8);
\draw [-, dashed, black] (-0.8+3.5,+0.8) -- (0.8+3.5,-0.8);

\node [align=center] at (-1+2.8,-1)
{$(K_4,\bold{-1})$};

\draw (-0.8+3.5+3.5,-0.8) [black,fill=black] circle (0.08 cm);
\draw (+0.8+3.5+3.5,-0.8) [black,fill=black] circle (0.08 cm);
\draw (-0.8+3.5+3.5,+0.8) [black,fill=black] circle (0.08 cm);
\draw (+0.8+3.5+3.5,+0.8) [black,fill=black] circle (0.08 cm);

\draw [-, black] (-0.8+3.5+3.5,-0.8) -- (0.8+3.5+3.5,-0.8);
\draw [-, black] (-0.8+3.5+3.5,-0.8) -- (0.8+3.5+3.5,0.8);
\draw [-, black] (-0.8+3.5+3.5,-0.8) -- (-0.8+3.5+3.5,+0.8);
\draw [-, black] (0.8+3.5+3.5,-0.8) -- (0.8+3.5+3.5,+0.8);

\draw [-,dashed, black] (-0.8+3.5+3.5,+0.8) -- (0.8+3.5+3.5,+0.8);
\draw [-, black] (-0.8+3.5+3.5,+0.8) -- (0.8+3.5+3.5,-0.8);

\node [align=center] at (-1+3+3.4,-1)
{$(K_4,\sigma_1)$};

\draw [-,dashed, black] (0+6.2+3.5,0-0.8) -- (1.6+6.2+3.5,1.6-0.8);
\draw [-, dashed,black] (0+6.2+3.5,0-0.8) -- (1.6+6.2+3.5,0-0.8);
\draw [-, dashed,black] (1.6+6.2+3.5,0-0.8) -- (1.6+6.2+3.5,1.6-0.8);
\draw [-, dashed,black] (0+6.2+3.5,0-0.8) -- (0+6.2+3.5,1.6-0.8);
\draw [-, dashed,black] (0+6.2+3.5,1.6-0.8) -- (1.6+6.2+3.5,1.6-0.8);

\draw (0+6.2+3.5,0-0.8) [black,fill=black] circle (0.08 cm);
\draw (1.6+6.2+3.5,0-0.8) [black,fill=black] circle (0.08 cm);
\draw (0+6.2+3.5,1.6-0.8) [black,fill=black] circle (0.08 cm);
\draw (1.6+6.2+3.5,1.6-0.8) [black,fill=black] circle (0.08 cm);

\node [align=center] at (-1+7+2.9,-1)
{$(D,\bold{-1})$};
\end{tikzpicture}
\end{center}\caption{The forbidden signed subgraphs for signed line graphs  (with $s=1$) whose underlying graph is a line graph.}\label{fig:fsegn}
\end{figure}

\begin{figure}[h]
\begin{center}
 \begin{tikzpicture}

\draw (-0.1+0,0) [black,fill=black] circle (0.08 cm);
\draw (-0.1-1,0.8) [black,fill=black] circle (0.08 cm);
\draw (-0.1+1,0.8) [black,fill=black] circle (0.08 cm);
\draw (-0.1+0,-0.8) [black,fill=black] circle (0.08 cm);

\draw [-,  black] (-0.1+0,0) -- (-0.1-1,0.8);
\draw [-,  black] (-0.1+0,0) -- (-0.1+1,0.8);
\draw [-,  black] (-0.1+0,0) -- (-0.1+0,-0.8);
\draw [-, black] (-0.1+1,0.8) -- (-0.1-1,0.8);

\node [align=center] at (-0.1-1+0.15,-1)
{$(P,\bold{1})$};

\draw (-0.8+3.5,-0.8) [black,fill=black] circle (0.08 cm);
\draw (+0.8+3.5,-0.8) [black,fill=black] circle (0.08 cm);
\draw (-0.8+3.5,+0.8) [black,fill=black] circle (0.08 cm);
\draw (+0.8+3.5,+0.8) [black,fill=black] circle (0.08 cm);

\draw [-,  black] (-0.8+3.5,-0.8) -- (0.8+3.5,-0.8);
\draw [-,  black] (-0.8+3.5,-0.8) -- (0.8+3.5,0.8);
\draw [-, black] (-0.8+3.5,-0.8) -- (-0.8+3.5,+0.8);
\draw [-,  black] (0.8+3.5,-0.8) -- (0.8+3.5,+0.8);

\draw [-,  black] (-0.8+3.5,+0.8) -- (0.8+3.5,+0.8);
\draw [-,  black] (-0.8+3.5,+0.8) -- (0.8+3.5,-0.8);

\node [align=center] at (-1+2.8,-1)
{$(K_4,\bold{1})$};

\draw (-0.8+3.5+3.5,-0.8) [black,fill=black] circle (0.08 cm);
\draw (+0.8+3.5+3.5,-0.8) [black,fill=black] circle (0.08 cm);
\draw (-0.8+3.5+3.5,+0.8) [black,fill=black] circle (0.08 cm);
\draw (+0.8+3.5+3.5,+0.8) [black,fill=black] circle (0.08 cm);

\draw [-, black] (-0.8+3.5+3.5,-0.8) -- (0.8+3.5+3.5,-0.8);
\draw [-, black] (-0.8+3.5+3.5,-0.8) -- (0.8+3.5+3.5,0.8);
\draw [-, black] (-0.8+3.5+3.5,-0.8) -- (-0.8+3.5+3.5,+0.8);
\draw [-, black] (0.8+3.5+3.5,-0.8) -- (0.8+3.5+3.5,+0.8);

\draw [-,dashed, black] (-0.8+3.5+3.5,+0.8) -- (0.8+3.5+3.5,+0.8);
\draw [-, black] (-0.8+3.5+3.5,+0.8) -- (0.8+3.5+3.5,-0.8);

\node [align=center] at (-1+3+3.4,-1)
{$(K_4,\sigma_1)$};

\draw [-, black] (0+6.2+3.5,0-0.8) -- (1.6+6.2+3.5,1.6-0.8);
\draw [-,black] (0+6.2+3.5,0-0.8) -- (1.6+6.2+3.5,0-0.8);
\draw [-,black] (1.6+6.2+3.5,0-0.8) -- (1.6+6.2+3.5,1.6-0.8);
\draw [-,black] (0+6.2+3.5,0-0.8) -- (0+6.2+3.5,1.6-0.8);
\draw [-,black] (0+6.2+3.5,1.6-0.8) -- (1.6+6.2+3.5,1.6-0.8);

\draw (0+6.2+3.5,0-0.8) [black,fill=black] circle (0.08 cm);
\draw (1.6+6.2+3.5,0-0.8) [black,fill=black] circle (0.08 cm);
\draw (0+6.2+3.5,1.6-0.8) [black,fill=black] circle (0.08 cm);
\draw (1.6+6.2+3.5,1.6-0.8) [black,fill=black] circle (0.08 cm);

\node [align=center] at (-1+7+2.9,-1)
{$(D,\bold{1})$};

\end{tikzpicture}
\end{center}\caption{The forbidden signed subgraphs for signed line graphs (with $s=-1$) whose underlying graph is a line graph.}\label{fig:fsegnmeno}
\end{figure}

From the classical theory we know that, if $L$ is a line graph, then $\Spec(L)\subseteq [-2,\infty)$.
From \cite{zasmat}  we know that if
$(L,\sigma)$ is a signed line graph with the choice  $s=-1$ then  $\Spec(L,\sigma)\subseteq (-\infty,2]$.
On the other hand, if $(L,\sigma)$ is a signed line graph with the choice $s=1$ then $\Spec(L,\sigma)\subseteq [-2,\infty)$ (see \cite{reff0,line} for more general results).
\\ Thanks to our characterization given in Theorem \ref{theo:1}, for a signed graph whose underlying graph is a line graph, it is possible to characterize the property of being a signed line graph by just looking at these spectral conditions.
\begin{theorem}\label{cor:spec}
Let $L=L(\Gamma)$ and let $\sigma$ be a signature of $L$. The following are equivalent.
\begin{enumerate}[(i)]
\item $(L,\sigma)$ is a signed line graph with the choice $s=1$.
\item $(L,\sigma)$ has no signed subgraph, induced by a vertex subset, which is switching isomorphic to one among $(P,\bold{-1})$, $(K_4,\bold{-1})$, $(K_4,\sigma_1)$, $(D,\bold{-1})$.
\item $\Spec(L,\sigma)\subseteq [-2,\infty)$.
\item $(L,-\sigma)$ is a signed line graph with the choice $s=-1$.
\item $(L,-\sigma)$ has no signed subgraph, induced by a vertex subset, which is switching isomorphic to one among $(P,\bold{1})$, $(K_4,\bold{1})$, $(K_4,\sigma_1)$, $(D,\bold{1})$.
\item $\Spec(L,-\sigma)\subseteq (-\infty,2]$.
\end{enumerate}
\end{theorem}
\begin{proof}
We will prove that conditions (i), (ii), (iii) are equivalent (the equivalence of (iv), (v), (vi) can be similarly proved). The proof will be concluded by noticing that (iii) is clearly equivalent to (vi).         \\
(i)$\implies$ (iii) It follows from \cite{reff0,line}.\\
(ii)$\implies$ (i) Suppose that (ii) holds and that $(L,\sigma)$ is not a signed line graph with $s=1$, so that $(L,\sigma)$ is not $\mathcal{Y}$-free by Theorem \ref{theo:1}. On the other hand, as $L$ is a line graph by the hypothesis, it must be $\mathcal{X}$-free: this implies that $(L,\sigma)$ is not $\mathcal{F}_1$-free and so it contains a signed subgraph, induced by some subset of $V_L$, which is switching isomorphic, by Remark \ref{rem:f}, to one among the graphs $(P,\bold{-1})$, $(K_4,\bold{-1})$, $(K_4,\sigma_1)$, $(D,\bold{-1})$ of Fig. \ref{fig:fsegn}. A contradiction.\\
(iii)$\implies$ (ii) Suppose $\Spec(L,\sigma)\subseteq [-2,\infty)$ and that $(L,\sigma)$ contains a signed subgraph, induced by a vertex subset, which is switching isomorphic to one among $(P,\bold{-1})$, $(K_4,\bold{-1})$, $(K_4,\sigma_1)$, $(D,\bold{-1})$.  An explicit computation shows that each of these $4$ signed graphs  has an eigenvalue less than $-2$. As a consequence of the Interlacing Theorem (see \cite[Theorem 1.3.11]{cve}) there is at least an eigenvalue of $(L,\sigma)$ with the same property, that is a contradiction.
\end{proof}

\begin{remark}\label{rem:vi}
Theorem \ref{cor:spec} can be compared with the results of \cite{meno2}, where the class of signed graphs represented by $D_\infty$ is characterized in terms of $49$ forbidden signed subgraphs $S_1,\ldots, S_{49}$. Notice that all such signed graphs have all eigenvalues greater than or equal to $-2$. The signed graphs in this class whose underlying graph is a line graph are exactly the  signed line graphs with $s=1$.
Now, one can remove from the list $S_1,\ldots, S_{49}$ those signed graphs whose underlying graph has subgraphs isomorphic to graphs in $\mathcal X$.
Those signed graphs in fact, by Beineke's characterization, never appear as subgraphs of a signed graph whose underlying graph is a line graph. It is easy to check that, after this deletion process, only the signed graphs $S_1,S_2,S_3,S_4$ remain: they are exactly switching-isomorphic copies of the signed graphs of Fig. \ref{fig:fsegn}.
\end{remark}
As a further consequence of Theorem \ref{cor:spec} and Corollary \ref{cor:cp}, we can deduce a characterization of (unsigned) line graphs with spectral radius at most $2$.
\begin{corollary}\label{cor:cp2}
A connected line graph $L$ with spectral radius at most $2$ must be a cycle or a path.
\end{corollary}
\begin{proof}
The maximal eigenvalue of a signed graph $(L,\sigma)$ must be less than or equal to that of the underlying graph $L$. It follows  that $\Spec(L,\sigma)\subseteq (-\infty,2]$ for every signature $\sigma$ of $L$. Then, by Theorem \ref{cor:spec},  the signed graph $(L,\sigma)$ is a signed line graph for every signature $\sigma$ of $L$.  By virtue of Corollary \ref{cor:cp}, the graph $L$ must be a cycle or a path.
\end{proof}
Notice that the same result can be deduced from \cite{smith}, where the graphs with spectral radius at most $2$ are classified.

\end{document}